\documentclass [a4paper,twoside,11pt]{article}

\usepackage[latin1]{inputenc}
\usepackage{amsmath,amsfonts,amssymb}
\usepackage{vmargin,graphicx,theorem}
\usepackage[english]{babel}
\usepackage{enumerate}

\setpapersize[portrait]{A4}
\setmarginsrb{2.2cm}{2cm}{2.2cm}{1cm}{0cm}{0.1cm}{0.5cm}{2cm}

\newcommand{\disp}{\displaystyle}

\newcommand{\dE}{\ensuremath{\mathbb{E}}}

\newcommand{\dL}{\ensuremath{\mathbb{L}}}

\newcommand{\dR}{\ensuremath{\mathbb{R}}}

\newtheorem{ethm}{Theorem}[section]

\newtheorem{ecor}[ethm]{Corollary}

\newtheorem{eprop}[ethm]{Proposition}

\newtheorem{elem}[ethm]{Lemma}

\newtheorem{edefi}[ethm]{Definition}

\newtheorem{erem}[ethm]{Remark}

\newtheorem{eex}[ethm]{Example}

\newcommand{\proofend}{~$\rhd$}
\newcommand{\proofbegin}{~$\lhd$}

\newenvironment{eproof}
               {\noindent {\emph{\textbf{Proof}}}\\\proofbegin~}
               {\proofend\\}

\newcommand{\p}[4]{{#3}\!\left#1{#4}\right#2}

\newcommand{\ABS}[1]{\ensuremath{{\left| #1 \right|}}} 
\newcommand{\PAR}[1]{\ensuremath{{\left(#1\right)}}} 
\newcommand{\BRA}[1]{\ensuremath{{\left\{#1\right\}}}} 
\newcommand{\NRM}[1]{\ensuremath{{\left\Vert #1\right\Vert}}} 

\renewcommand{\phi}{\varphi}

\renewcommand{\geq}{\geqslant}

\newcommand{\varf}[1]{\mathbf{Var}_{#1}}
\newcommand{\entf}[1]{\mathbf{Ent}_{#1}}

\newcommand{\ent}[2]{\p(){\entf{#1}}{#2}}

\newcommand{\var}[2]{\p(){\varf{#1}}{#2}}

\def\disp{\displaystyle}

\newcommand{\R}{\dR}

\newcommand{\Pt}[1][t]{\ensuremath{\mathbf{P_{\!#1}}}}
\newcommand{\Ps}[1][s]{\ensuremath{\mathbf{P_{\!#1}}}}

\renewcommand{\P}{\mathbf{P}}

\newcommand{\al}{\alpha}

\newcommand{\e}{\varepsilon}

\newcommand{\1}{\hbox{1}\!\!\hbox{I}}

\newcommand{\osc}[1]{{\bf Osc}\PAR{#1}} 

\begin{document}

\title{\sl Weak logarithmic Sobolev inequalities and entropic convergence}
\author{P. Cattiaux,  I. Gentil and  A. Guillin }

\date{\today}
\maketitle\thispagestyle{empty}

\begin{abstract}
In this paper we introduce and study a weakened form of logarithmic Sobolev inequalities in
connection with various others functional inequalities (weak Poincar\'e inequalities, general Beckner inequalities...).  We also discuss the quantitative behaviour of relative entropy along a symmetric diffusion semi-group. In particular, we exhibit an example where Poincar\'e inequality can not be used for deriving entropic convergence whence weak logarithmic Sobolev inequality ensures the result.
\end{abstract}

\bigskip


\par\vspace{40pt}
\section{Introduction}\label{sec-intro}

Since the beginning of the nineties, functional inequalities (Poincar\'e, logarithmic (or F-)
Sobolev, Beckner's like, transportation) turned to be a powerful tool for studying various problems
in Probability theory and in Statistics: uniform ergodic theory, concentration of measure,
empirical processes, statistical mechanics, particle systems for non linear p.d.e.'s, stochastic
analysis on path spaces, rate of convergence of p.d.e....

Among such functional inequalities, Poincar\'e inequality and its generalizations (weak and super
Poincar\'e) deserved particular interest, as they are the most efficient tool for the study of
isoperimetry, concentration of measure and $\dL^2$ long time behavior (see e.g.
\cite{r-w,w1,w2,ca-ba-ro2,ca-ba-ro3}). However (except the usual Poincar\'e inequality) they are
not easily tensorizable nor perturbation stable. That is why super-Poincar\'e inequalities have to
be compared with (generalized) Beckner's inequalities or with additive $\phi$-Sobolev inequalities
(see \cite{w2,ca-ba-ro3,ca-ba-ro}).

But for some aspects, generalized Poincar\'e inequalities are insufficient. Indeed $\dL^2$ controls
are not well suited in various situations (statistical mechanics, non linear p.d.e), where entropic
controls are more natural. It is thus interesting to look at generalizations of Gross
logarithmic Sobolev inequality. In this paper we shall investigate weak logarithmic Sobolev inequalities
(the ``super'' logarithmic Sobolev inequalities have already been investigated by Davies and Simon, or
R\"{o}ckner and Wang).

In order to better understand the previous introduction and what can be expected, let us introduce
some definitions and recall some known facts. In all the paper $M$ denotes a Riemannian manifold and $\mu$ denotes an absolutely
continuous probability measure with respect to the surface measure.  We also assume that $\mu$ is
symmetric for a diffusion semi-group $P_t$ associated to a non explosive
diffusion process.

Let $H^1(M, \mu)$ be the closure of $C_b^\infty(M)$ (the space of
infinitely differentiable functions $f$ on $M$ with all $|\nabla^nf|, n\ge 0$ bounded) w.r.t. the norm
$
\sqrt{\mu(|f|^2 + |\nabla f|^2)}.
$

\begin{edefi}
\label{def-wp} We say that the measure $\mu$ satisfies  a weak Poincar\'e inequality, {\bf WPI}, if
there exists a non-increasing function $\beta_{WP} : (0,+\infty)\rightarrow\dR^+$, such that for
all $s>0$ and any bounded function $f\in H^1(M,\mu)$,
\begin{equation}
\tag{{\bf WPI}} \varf{\mu}(f):=\int f^2d\mu-\PAR{\int fd\mu}^2 \leq\beta_{WP}(s)\int\ABS{\nabla f}^2d\mu+s \, {\bf Osc}^2(f),
\end{equation}
where $\osc{f}=\sup f-\inf f$.
 \end{edefi}

Weak Poincar\'e inequalities have been introduced by R\"{o}ckner and Wang in \cite{r-w}. If
$\beta_{WP}$ is bounded, we recover the (classical) Poincar\'e inequality, while if
$\beta_{WP}(s)\to\infty$ as $s\to 0$ we obtain a weaker inequality.

Actually, as shown in \cite{r-w} any Boltzman measure ($d\mu = e^{-V} \, dx$) on $\dR^n$ with a
locally bounded potential $V$ satisfies some {\bf WPI} (the result extends to any manifold with
Ricci curvature bounded from below by a possibly negative constant, according to Theorem 3.1 in
\cite{r-w} and the local Poincar\'e inequality shown by Buser \cite{Bus} in this framework). {\bf
WPI} furnishes an isoperimetric inequality, hence (sub-exponential) concentration of measure (see
\cite{r-w,ca-ba-ro2}). It also allows to describe non exponential decay of the $\dL^2$ norm of the
semi group, i.e. {\bf WPI} is linked to inequalities like
$$
\forall t\geq0,\quad\varf{\mu}(P_t f) \, \leq \, \xi(t) \, {\bf Osc}^2(f),
$$
 for some adapting function $\xi$
(relations between $\beta_{WP}$ and $\xi$ will be recalled later). Recall that a uniform decay of
the Variance, is equivalent to its exponential decay which is equivalent to the usual Poincar\'e
inequality. Let us note that a multiplicative form of weak Poincar\'e inequality (namely
$\beta(s)=s^{2/p}$ and choose $s$ such that each term of the right hand side is of the same order)
appears first in works of Liggett \cite{Lig} to prove an algebraic convergence in $L^2$ of some
spin system dynamic.

If we replace the variance by the entropy
 the latter argument is still true. Indeed (at least for bounded below curvature) an
uniform decay of $\ent{\mu}{P_t h}$ is equivalent to its exponential decay which is equivalent to
the logarithmic Sobolev inequality. In order to describe non exponential decays, it is thus natural to
introduce the following definition:

\begin{edefi}
\label{def-wls} We say that the measure $\mu$ satisfies  a weak logarithmic Sobolev inequality, {\bf
WLSI}, if there exists a non-increasing function $\beta_{WL} : (0,+\infty)\rightarrow\dR^+$, such
that for all $s>0$ and any bounded function $f\in H^1(M,\mu)$,
\begin{equation}
\tag{{\bf WLSI}} \ent{\mu}{f^2}:=\int f^2 \log\PAR{\frac{f^2}{\int f^2d\mu}}d\mu
\leq\beta_{WL}(s)\int\ABS{\nabla f}^2d\mu+s \, {\bf Osc}^2(f) \, .
\end{equation}
\end{edefi}

Remark that {\bf WPI} is translation invariant. Hence it is enough to check it for non negative
functions $f$ and  for such functions we get  $\varf{\mu}(f)\leq \ent{\mu}{f^2}$. Hence {\bf WLSI} is stronger
than {\bf WPI} (we shall prove a more interesting result), and we can expect that {\bf WLSI} (with
a non bounded $\beta_{WL}$) allows to describe all the sub-gaussian measures, in particular all
super-exponential (and sub-gaussian measures) for which a strong form of Poincar\'e inequality
holds. Remark that, as for weak Poincar\'e inequalities, multiplicative forms of the weak logarithmic Sobolev inequality appears first under the name of log-Nash inequality to study the decay of semigroup in the case of Gibbs measures, see Bertini-Zegarlinski \cite{BZ1,BZ2} or Zegarlinski \cite{Z}.

\begin{erem}\label{rema}
One can easily check that $\varf{\mu}(f)\le {1\over 4}{\bf Osc}^2(f)$ so that we may assume that $\beta_{WP}(s)=0$ as soon as $s\ge 1/4$. In fact, one can also prove  $\ent{\mu}{f^2}\le {1\over e}{\bf Osc}^2(f)$ and thus we can suppose that $\beta_{WL}(s)=0$ for $s\ge 1/e$.\\
Hence for {\bf WPI} and {\bf WLSI} what is important is the behaviour of $\beta$ near 0.
\end{erem}

In order to understand the picture and to compare all these inequalities we shall call upon another class of inequalities, namely measure-capacity inequalities introduced by Maz'ya \cite{Maz}. Then these inequalities are  extensively used in this context  \cite{barthe-roberto,chen05,ca-ba-ro,ca-ba-ro2,ca-ba-ro3}. Given measurable sets $A\subset\Omega$ the capacity $Cap_\mu(A,\Omega)$, is defined as follow:
$$
Cap_\mu(A,\Omega):=\inf\BRA{\int\ABS{\nabla f}^2d\mu;\,\1_A\leq f\leq\1_\Omega},
$$
where the infimum is taken over all function $f\in H^1(M,\mu)$. By convention, if the set of function $f\in  H^1(M,\mu)$ such that $\1_A\leq f\leq\1_\Omega$ is empty then we note $Cap_\mu(A,\Omega)=+\infty$.   We refer to Maz'ya \cite{Maz} and Grigor'yan \cite{Gri} for further discussion on capacities. The capacity defined by Maz'ya seems to be  a little different but they are similar. If now $A$ satisfies $\mu(A)\leq 1/2$
we note
\begin{equation}
\label{def-cap}
Cap_\mu(A):=\inf\BRA{Cap_\mu(A,\Omega);\,A\subset\Omega,\,\mu(\Omega)\leq1/2}.
\end{equation}
A measure-capacity inequality is an inequality of the form
\begin{equation}
\label{def-mc}
\frac{\mu(A)}{\gamma(\mu(A))}\leq Cap_\mu(A),
\end{equation}
for some function $\gamma$. They are in a sense universal, since
they only involve the energy (Dirichlet form) and the measure. Furthermore, a remarkable feature is
that most of known inequalities involving various functionals (variance, $p$-variance, $F$
functions of $F$-Sobolev inequalities, entropy etc...) can be compared (in a non sharp form) with
some measure-capacity inequalities.

We shall thus start by characterizing {\bf WLSI} via measure-capacity inequalities. Then we will
study the one dimensional case, in the spirit of Muckenhoupt or Bobkov-G\"{o}tze criteria for
Poincar\'e or logarithmic Sobolev inequalities (see e.g. \cite{logsob} chapter 6). We shall then discuss in details the
relationship between {\bf WLSI} and the generalized Poincar\'e inequalities. Finally we shall
discuss various properties and consequences of {\bf WLSI}. In the final sections, we
study in details the decay of entropy for large time. In particular we show that for a $\mu$
reversible gradient diffusion process, very mild conditions on the initial law are sufficient to
ensure an entropic decay of type $e^{-t^\beta}$ when $\mu$ satisfies interpolating inequalities between
Poincar\'e and Gross introduced by Latala and Oleszkiewicz \cite{LO00}, those conditions preventing estimation via Poincar\'e inequalities. We also give the elements
to compute this decay under general {\bf WLSI}. The particular case of the double sided exponential measure is detailed.
\medskip

Let us finally remark that the limitation to finite dimensional space is only instrumental and the main results would be readily extendable to infinite dimensional space with capacity defined to suitable Dirichlet forms (assuming for example the existence of a {\it carr\'e du champ} operator).

\section{Weak logarithmic Sobolev inequalities}
\label{sec-wls}

\subsection{Characterization via capacity-measure condition}

We start this section by characterizing {\bf WLSI} in terms of measure-capacity inequalities.
\begin{ethm}
\label{thm-cap} Assume that the measure $\mu$ satisfies a {\bf WLSI}  with
function $\beta_{WL}$, then for every $A\subset M$ such that $\mu(A)\leq 1/2$,
$$
\forall s>0, \quad\frac{\mu(A)\log\PAR{1+\frac{1}{2\mu(A)}}-s}{\beta_{WL}(s)}\leq Cap_\mu(A).
$$
\end{ethm}

\begin{eproof}
Let $A\subset\Omega$ with $\mu(\Omega)\leq1/2$ and let $f$ be a locally Lipschitz function
satisfying $\1_A\leq f\le \1_\Omega$.  The variational definition of the entropy implies
$$
\ent{\mu}{f^2}\geq \int f^2 gd\mu,
$$
for all $g$ such that $\int e^gd\mu \leq1$. Apply this inequality with
\begin{equation*}
g=\left\{
\begin{array}{ll}
\log\PAR{1+\frac{1}{2\mu(A)}}& \text{ on } A\\
0 & \text{ on }\Omega\backslash A\\
-\infty&\text{ on }\Omega^c
 \end{array}
\right.
\end{equation*}
which satisfies $\int e^gd\mu\leq1$. It yields
$\ent{\mu}{f^2}\geq\mu(A)\log\left(1+\frac{1}{2\mu(A)}\right)$.

Therefore by the weak logarithmic Sobolev inequality and the definition of the capacity  we obtain
$$
\mu(A)\log\left(1+\frac{1}{2\mu(A)}\right)\leq \beta_{WL}(s)Cap_\mu(A,\Omega)+s.
$$
Taking the infimum over sets $\Omega$ with measure at most 1/2 and containing $A$ we obtain
$$
\forall s>0,\quad\frac{\mu(A)\log\PAR{1+\frac{1}{2\mu(A)}}-s}{\beta_{WL}(s)}\leq Cap_\mu(A).
$$
\end{eproof}

\begin{ethm}
\label{thm-wls} Let $\beta:(0,+\infty)\rightarrow\dR^+$ be non-increasing function such that for every
$A\subset M$ with $\mu(A)\leq 1/2$ one has
\begin{equation}
\label{eq-wls} \forall
s>0,\,\,\,\,\,\,\frac{\mu(A)\log\PAR{1+\frac{e^2}{\mu(A)}}-s}{\beta(s)}\leq Cap_\mu(A).
\end{equation}

Then the measure $\mu$ satisfies a {\bf WLSI}  with the function $\beta_{WL}(s)=16 \beta(3s/14)$, for $s>0$.
\end{ethm}

\begin{eproof}
  Let a bounded function $f\in H^1(M,\mu)$, we will prove that
\begin{equation}
\label{eq-thm2}
\forall s>0,\quad \ent{\mu}{f^2}\leq 16\beta(s)\int\ABS{\nabla f}^2d\mu+ 14s/3 \, {\bf Osc}^2(f).
\end{equation}

Let $m$ be a median of $f$ under $\mu$ and let $\Omega_+=\BRA{f>m}$,
$\Omega_-=\BRA{f<m}$. Then, using the argument of Lemma 5 in \cite{barthe-roberto}, we obtain
\begin{multline}
\label{eq-bbr}
\ent{\mu}{f^2}\leq\sup\BRA{\int F_+^2hd\mu;\,\,h\geq0,\,\,\int e^hd\mu\leq e^2+1}\\
+\sup\BRA{\int F_-^2hd\mu;\,\,h\geq0,\,\,\int e^hd\mu\leq e^2+1},
\end{multline}
where $F_+=(f-m)\1_{\Omega_+}$ and $F_-=(f-m)\1_{\Omega_-}$.

We will study the first term in the right hand side, the second one will be treated by the
same method.

There are two cases depending on the value of $s$. Let $s_1:=\frac{1}{2}\log\PAR{1+2e^2}$, and assume that
$s\in(0,s_1)$. Let define $c$ by
$$
c=\inf\BRA{t\geq 0,\quad \mu(F_+^2>t)\log\PAR{1+\frac{e^2}{\mu(F_+^2>t)}}\leq s}.
$$
If $c=0$ then one get that for some constant $C$
$$
\sup\BRA{\int F_+^2hd\mu;\,\,h\geq0,\,\,\int e^hd\mu\leq e^2+1}\leq s  \log(1+e^2)\NRM{F_+}_\infty^2
$$
and the problem is solved on that case. If now $c>0$, since $\mu$ is absolutely continuous with respect to the surface measure of the Riemannian manifold $M$, one can find $\Omega_0$ such that
$\BRA{F^2_+>c}\subset  \Omega_0\subset \BRA{F^2_+\geq c}$ and
\begin{equation}
\label{def-c}
\mu(\Omega_0)\log\PAR{1+\frac{e^2}{\mu(\Omega_0)}}= s.
\end{equation}
Note that the function $x \mapsto x\log(1+e^2/x)$ is increasing on $(0,\infty)$, and realize a bijection between $(0,1/2]$
and $(0,s_1]$.

Pick some $\rho\in(0,1)$ and introduce for any $ k>0$, $\Omega_k=\BRA{F_+^2\geq c \rho^k}$.
The sequence $\PAR{\Omega_k}_k$ is increasing so that, for every function $h\geq0$,
$$
\int F_+^2hd\mu= \int_{\Omega_0} F_+^2hd\mu+\sum_{k>0}\int_{\Omega_k\backslash\Omega_{k-1}} F_+^2hd\mu.
$$
For the first term we get
$$
\int_{\Omega_0} F_+^2hd\mu\leq {\bf Osc}^2\PAR{f}\int_{\Omega_0} h \, d\mu,
$$
then Lemma 6 of \cite{barthe-roberto} implies that
$$
\sup\BRA{\int_{\Omega_0} hd\mu;\,\,h
\geq0,\,\,\int e^hd\mu\leq e^2+1}=\mu(\Omega_0)\log\PAR{1+\frac{e^2}{\mu(\Omega_0)}}.
$$
So that, using the definition of $c$ (equality~\eqref{def-c}) we get
$$
\sup\BRA{\int_{\Omega_0} F_+^2hd\mu;\,\,h\geq0,\,\,\int e^hd\mu\leq e^2+1}\leq s \, {\bf Osc}^2(f).
$$

For the second term we have for all $k>0$, due to the fact that  $c\rho^{k}\leq F_+^2\leq c\rho^{k-1}$ on
$\Omega_k\backslash\Omega_{k-1}$,
$$
\int_{\Omega_k\backslash\Omega_{k-1}} F_+^2hd\mu\leq c \rho^{k-1}
\int_{\Omega_k\backslash\Omega_{k-1}} hd\mu.
$$
Then we obtain using again Lemma 6 of \cite{barthe-roberto}, for any $k>0$,
\begin{equation*}
\sup\BRA{\int_{\Omega_k\backslash\Omega_{k-1}} F_+^2hd\mu;\,\,h\geq0,\,\,\int e^hd\mu\leq e^2+1}
\leq c \rho^{k-1}\mu(\Omega_k\backslash\Omega_{k-1})\log\PAR{1+\frac{e^2}{\mu(\Omega_k\backslash\Omega_{k-1})}}.
\end{equation*}

Using now inequality~\eqref{eq-wls} we get
$$
c\rho^{k-1}\PAR{\mu(\Omega_k\backslash\Omega_{k-1})\log\PAR{1+\frac{e^2}{\mu(\Omega_k\backslash\Omega_{k-1})}}}
\leq c\rho^{k-1}\beta(s)Cap_\mu(\Omega_k\backslash\Omega_{k-1})+sc\rho^{k-1}.
$$

Let set for any $k>0$,
$$
g_k=\min\BRA{1,\PAR{\frac{F_+-\sqrt{c\rho^{k+1}}}{\sqrt{c\rho^{k}} - \sqrt{c\rho^{k+1}}}}_+},
$$
so that we have $\1_{\Omega_{k}}\leq g_k\leq \1_{\Omega_{k+1}}$ with $\mu(\Omega_+)\le 1/2$. This implies, using the
 definition of $Cap_\mu(\Omega_k\backslash\Omega_{k-1})$ (see~\eqref{def-cap}),
$$
c \rho^{k-1} Cap_\mu(\Omega_k\backslash\Omega_{k-1})     \leq
\frac{1}{\rho(1-\sqrt{\rho})^2} \int_{\Omega_{k+1}\backslash\Omega_{k}}\ABS{\nabla F_+}^2d\mu.
$$

\medskip

Note that the constant $c$ satisfies $c\leq\NRM{F_+}_\infty^2\leq {\bf Osc}^2\PAR{f}$. We can now
finish the proof in the case $s\in(0,s_1)$,
\begin{eqnarray*}
\disp\sup\BRA{\int F_+^2hd\mu;\,\,h\geq0,\,\,\int e^hd\mu\leq e^2+1}\!\!\!\! & \leq  &
\sup\BRA{\int_{\Omega_0} F_+^2hd\mu;\,\,h\geq0,\,\,\int e^hd\mu}+\\
&&\qquad \disp \sum_{k>0}\sup\BRA{\int_{\Omega_{k+1}\backslash\Omega_{k}} F_+^2hd\mu;\,\,h\geq0,\,\,\int e^hd\mu}\\
&\leq & s \, \disp{\bf Osc}^2(f)+\sum_{k>0}sc\rho^{k-1} +\\
& &\qquad \sum_{k>0}\frac{1}{\rho(1-\sqrt{\rho})^2} \int_{\Omega_{k+1}\backslash\Omega_{k}}\ABS{\nabla F_+}^2d\mu\\
&\leq &\disp \frac{\beta(s)}{\rho(1-\sqrt{\rho})^2}\int\ABS{\nabla F_+}^2d\mu+
s\frac{2-\rho}{1-\rho}{\bf Osc}^2\PAR{f}.
\end{eqnarray*}

Using inequality~\eqref{eq-bbr} and the previous inequality for $F_-$ we get
\begin{equation}
\label{eq-pr}
\forall s\in(0,s_1),\quad
\ent{\mu}{f^2}\leq\frac{\beta(s)}{\rho(1-\sqrt{\rho})^2}\int\ABS{\nabla f}^2d\mu+ 2s\frac{2-\rho}{1-\rho}{\bf
Osc}^2\PAR{f},
\end{equation}
for all $\rho\in(0,1)$. Choosing $\rho=1/4$ furnishes inequality~\eqref{eq-thm2} for any $s\in(0,s_1)$.

Assume now that $s\geq s_1$, then take $c=0$ and  we get
\begin{equation*}
\mu(\Omega_0)\log\PAR{1+\frac{e^2}{\mu(\Omega_0)}}\leq s,
\end{equation*}
and the same argument used for $s\in(0,s_1)$  implies
\begin{equation}
\label{eq-pr2}
\forall s\geq s_1,\quad
\ent{\mu}{f^2}\leq 2s{\bf Osc}^2\PAR{f}.
\end{equation}
Then inequality~\eqref{eq-pr2} and the previous result implies inequality~\eqref{eq-thm2} for any $s>0$.

Note that we do not obtain the optimal function $\beta_{WL}(s)$ for $s$ large, but, as explained in remark~\ref{rema},
this is not important for the {\bf WLSI}.
\end{eproof}

\begin{erem}
\label{rem-eq}
The following two inequalities hold
$$
\frac{\frac{\mu(A)}{2}\log\PAR{1+\frac{1}{2\mu(A)}}}{\beta_{WL}\PAR{\frac{\mu(A)}{2}\log\PAR{1+\frac{1}{2\mu(A)}}}}\leq
\sup_{s>0} \BRA{\frac{\mu(A)\log\PAR{1+\frac{1}{2\mu(A)}}-s}{\beta_{WL}(s)}}\leq
\frac{{\mu(A)}\log\PAR{1+\frac{1}{2\mu(A)}}}{\beta_{WL}\PAR{{\mu(A)}\log\PAR{1+\frac{1}{2\mu(A)}}}
}
$$
and
\begin{multline}
\label{eq-ii}
\frac{{\frac{\mu(A)}{2}\log\PAR{1+\frac{e^2}{\mu(A)}}}}
{\beta_{WL}\PAR{{\frac{\mu(A)}{2}\log\PAR{1+\frac{e^2}{\mu(A)}}}}}\leq\\
\sup_{s>0}\BRA{\frac{\mu(A)\log\PAR{1+\frac{e^2}{\mu(A)}}-s}{\beta_{WL}(s)}}\leq
\frac{{\mu(A)\log\PAR{1+\frac{e^2}{\mu(A)}}}}{\beta_{WL}\PAR{{\mu(A)\log\PAR{1+\frac{e^2}{\mu(A)}}}}}.
\end{multline}
Proofs of these inequalities are the same as in~\cite[Theorem 1]{ca-ba-ro2}. The lower bounds
of these inequalities
correspond to a specific choice, $s=\frac{\mu(A)}{2}\log\PAR{1+\frac{1}{2\mu(A)}}$ for the first one and
$s=\frac{\mu(A)}{2}\log\PAR{1+\frac{e^2}{\mu(A)}}$ for the second one.
For the upper bound of the first inequality we use the fact that
$$
\sup_{s>0}\BRA{\frac{\mu(A)\log\PAR{1+\frac{1}{2\mu(A)}}-s}{\beta_{WL}(s)}}\leq
\sup_{0<s<\mu(A)\log\PAR{1+\frac{1}{2\mu(A)}}}\BRA{\frac{\mu(A)\log\PAR{1+\frac{1}{2\mu(A)}}}{\beta_{WL}(s)}},
$$
and the  non-increasing property of $\beta$ gives the result. The method holds for the second inequality.
\end{erem}

\subsection{A Hardy like criterion on $\dR$}

\begin{eprop}
\label{thm-dim1}
   Let $\mu$ be a probability measure on $\dR$. Assume that $\mu$ is absolutely
continuous with respect to Lebesgue measure and denote by $\rho_\mu$ its density. Let $m$ be a
median of $\mu$ and $\beta_{WL}:(0,\infty)\to \dR^{+}$  be non-increasing. Let $C$
be the optimal constant such that for all $f\in H^1(\dR,\mu)$,
$$
\forall s>0,\quad
\ent{\mu}{f^2} \le C\beta_{WL}(s) \int |\nabla f|^{2}d\mu +s \, {\bf Osc}^2(f).
$$
Then we get  $\max(b_-,b_+)  \le  C \le \max (B_-,B_+)$, where
\begin{eqnarray}
\nonumber
   b_+&:=&\sup_{x>m} \frac{\frac{\mu([x,+\infty))}{2}\log\PAR{1+\frac{1}{2\mu([x,+\infty))}}}
{{\beta_{WL}\PAR{\frac{\mu([x,+\infty))}{2}\log\PAR{1+\frac{1}{2\mu([x,+\infty))}}}}}\int_m^{x}\frac{1}{\rho_\mu} \\
\nonumber
  b_-&:=&\sup_{x<m} \frac{\frac{\mu((-\infty,x])}{2}\log\PAR{1+\frac{1}{2\mu((-\infty,x])}}}
{{\beta_{WL}\PAR{\frac{\mu((-\infty,x])}{2}\log\PAR{1+\frac{1}{2\mu((-\infty,x])}}}}}
              \int_x^{m}\frac{1}{\rho_\mu} \\
\label{eq-B+}
    B_+&:=&\sup_{x>m}  \frac{16{\mu([x,+\infty))}\log\PAR{1+\frac{e^2}{\mu([x,+\infty))}}}
{{\beta_{WL}\PAR{{14\over 3}{\mu([x,+\infty))}\log\PAR{1+\frac{e^2}{\mu([x,+\infty))}}}}}\int_m^{x}\frac{1}{\rho_\mu} \\
\nonumber
        B_-&:=&\sup_{x<m} \frac{16{\mu((-\infty,x])}\log\PAR{1+\frac{e^2}{\mu((-\infty,x])}}}
{{\beta_{WL}\PAR{ {14\over 3}{\mu((-\infty,x])}\log\PAR{1+\frac{e^2}{\mu((-\infty,x])}}}}}
              \int_x^{m}\frac{1}{\rho_\mu}
\end{eqnarray}
\end{eprop}

\begin{eproof}
The proof of the lower bound on $C$ is exactly the same as in \cite[Theorem~3]{ca-ba-ro2} using
Theorem~\ref{thm-cap} and  Remark~\ref{rem-eq}.

For the upper bound  denote $F_+=(f-f(m))\1_{[m,+\infty)}$ and $F_-=(f-f(m))\1_{(-\infty,m]}$. Then
$$
\ent{\mu}{f^2}\leq\ent{\mu}{F_+^2}+\ent{\mu}{F_-^2}.
$$
We work separately with the two terms and explain the arguments for $\ent{\mu}{F_+^2}$ only. We
follow the method of proof in \cite[Theorem~3]{ca-ba-ro2}.

Using equality~\eqref{eq-B+} we get
$$
\forall x>m,\quad\frac{16{\mu([x,+\infty))}\log\PAR{1+\frac{e^2}{\mu([x,+\infty))}}}
{{\beta_{WL}\PAR{{14\over 3}{\mu([x,+\infty))}\log\PAR{1+\frac{e^2}{\mu([x,+\infty))}}}}}\int_m^{x}\frac{1}{\rho_\mu}\leq
B_+.
$$
This means that
$$
\forall x>m,\quad\frac{16{\mu([x,+\infty))}\log\PAR{1+\frac{e^2}{\mu([x,+\infty))}}}
{{B_+\beta_{WL}\PAR{{14\over 3}{\mu([x,+\infty))}\log\PAR{1+\frac{e^2}{\mu([x,+\infty))}}}}}\leq
Cap_\mu([x,+\infty),[m,+\infty)).
$$
If $A\subset [m,+\infty)$ then $Cap_\mu(A,[m,+\infty))=Cap_\mu([\inf
A,+\infty),[m,+\infty))$ (see for example \cite[Sec. 4]{barthe-roberto}). The  function
$$
t\mapsto \frac{16t\log\PAR{1+\frac{e^2}{t}}}{\beta_{WL}\PAR{{14\over3}t\log\PAR{1+\frac{e^2}{t}}}}
$$
is increasing on $(0,\infty)$, so we get
$$
\forall A\subset [m,+\infty),\quad
\frac{16{\mu(A)}\log\PAR{1+\frac{e^2}{\mu(A)}}} {{B_+ \,
\beta_{WL}\PAR{{14\over 3}{\mu(A)}\log\PAR{1+\frac{e^2}{\mu(A)}}}}}\leq Cap_\mu(A,[m,+\infty)).
$$
Using now inequality~\eqref{eq-ii} one has for all $A\subset [m,+\infty)$,
$$
\sup_{s>0}\BRA{16\frac{\mu(A)\log\PAR{1+\frac{e^2}{\mu(A)}}-s}{B_+ \, \beta_{WL}({14\over 3}s)}}\leq
Cap_\mu(A,[m,+\infty)),
$$
and then by the same argument as in Theorem~\ref{thm-wls} one has
\begin{equation*}
\ent{\mu}{F_+^2}\leq B_+\beta_{WL}(s)\int\ABS{\nabla F_+}^2d\mu+
s\osc{f}^2.
\end{equation*}
It follows that $C\leq B_+$. The same argument gives also $C\leq B_-$ and the proposition is
proved.
\end{eproof}

\begin{ecor}
\label{cor_wls}
  Let $\Phi$ be a function on $\dR$ such that $d\mu_\Phi(x):=e^{-\Phi(x)}dx,$ $x\in \dR$ is a probability measure and
 let $\varepsilon \in (0,1)$.

Assume that
there exists an interval $I=(x_0,x_1)$ containing a median $m$ of $\mu$ such
that $|\Phi|$ is bounded on $I$, and $\Phi$ is twice differentiable outside $I$
with for any $x\not\in I$,
\begin{eqnarray}
\nonumber
&&\disp\Phi'(x)\neq 0,\disp \frac{|\Phi''(x)|}{\Phi'(x)^{2}}\le 1-\varepsilon \text{ and }\\
\label{eq-ap}
&&\disp A'\Phi({x})\leq \Phi(x)+\log\ABS{\Phi'(x)}\leq A\Phi(x),
\end{eqnarray}
for  some constants $A,A'>0$.

Let  $\beta$ be a non-increasing function on $(0,\infty)$.
 Assume that there exists $c>0$ such that
for all $x\not\in I$ it holds
$$
\frac{\Phi(x)}{\Phi'(x)^2}\leq c \beta \PAR{\frac{A e^{-\Phi(x)}\Phi(x)}{\ABS{\Phi'(x)}}}.
$$
Then $\mu_\Phi$ satisfies a {\bf WLSI} with function $C\beta$ for some
constant $C>0$.
\end{ecor}

\begin{eproof}
Corollary~2.4 of \cite{ca-ba-ro2} gives for $x\ge x_1$,
$$ \mu([x,+\infty))\le \frac{e^{-\Phi(x)}}{\varepsilon \Phi'(x)}
\le \frac{2-\varepsilon}{\varepsilon}\mu([x,+\infty)).$$
Then using Proposition~\ref{thm-dim1} and inequality~\eqref{eq-ap} we obtain the result.
\end{eproof}

\begin{eex}
\label{ex1}
Let us give two examples:
\begin{itemize}
\item For $\alpha>0$, the measure $dm_\alpha(t)= \alpha(1+|t|)^{-1-\alpha}dt/2,\, t\in \dR$
satisfies the {\bf WLSI} with the function
$$
\forall s>0,\quad\beta_{WL}(s)=C\frac{\PAR{\log 1/s}^{1+2/\al}}{s^{2/\al}},
$$
for some constant $C>0$.
\item Let  $\al\in(0,2)$ and defined  the probability measure $d\mu_\al(t)=Z_\al
e^{-|t|^{\al}}dt,\; t\in \dR$, ($Z_\al$ is a normalization constant). Then $\mu_\al$  satisfies the {\bf WLSI}
with the function
$$
\forall s>0,\quad\beta_{WL}(s)=C\PAR{\log 1/s}^{(2-\al)/\al},
$$
for some  $C>0$.

 Contrary to the {\bf WPI}, one can
study the case $\al\in(1,2]$. In particular for $\al=2$ we get that $\beta_{WL}$ is bounded, i.e.
we recover (with a non sharp constant) the classical logarithmic Sobolev inequality for the gaussian
measure.
\end{itemize}
\end{eex}

\section{Weak Logarithmic Sobolev inequalities and generalized Poincar\'e inequalities}
\label{sec-poinc}

\subsection{Link with weak Poincar\'e inequalities and classical Poincar\'e inequality }


Barthe, Cattiaux and Roberto  investigated in \cite{ca-ba-ro2} the
 measure-capacity criterion for {\bf WPI}.
Their results read as follows: {\bf WPI} with a function $\beta_{WP}$ implies a measure-capacity
inequality with $\gamma(u)=4 \beta_{WP}(u/4)$ (see inequality~\eqref{def-mc}) while a measure-capacity inequality with
non-increasing function $\gamma$ implies {\bf WPI} with $\beta_{WP}=12 \gamma$ (we may assume that
$\gamma(u)=\gamma(1/2)$ for $u\geq 1/2$). Comparing with Theorem \ref{thm-cap} and Theorem
\ref{thm-wls}, we can state:

 \begin{eprop}
\label{prop-dd}
Assume that a probability measure $\mu$ satisfies a {\bf WLSI}
with function $\beta_{WL}$ then $\mu$ satisfies a {\bf WPI} with function
 $\beta_{WP}$
defined by
\begin{equation}
\label{eq-dd}
\forall s>0,\quad\beta_{WP}(s)={24\beta_{WL}\PAR{\frac{s}{2}\log\left(1+{1\over 2s}\right)}\over
\log\left(1+{1\over 2s}\right)}.
\end{equation}

 Conversely, a {\bf WPI} with  function $\beta_{WP}$ implies a {\bf WLSI}
 with function $\beta_{WL}$, defined by,
\begin{eqnarray}
\label{eq-ddd}
\left\{
\begin{array}{l}
\disp\forall s\in(0,s_0),\quad\disp\beta_{WL}(s)=c'\beta_{WP}\PAR{c\frac{s}{\log\PAR{1/s}}}\log\PAR{1/s},\\
\disp\forall s\geq s_0,\quad\beta_{WL}(s)=c'\beta_{WP}\PAR{c\frac{s_0}{\log\PAR{1/s_0}}}\log\PAR{1/s_0},
\end{array}
\right.
\end{eqnarray}
for some universal constants $c,c',s_0>0$.

Finally assume that $\mu$ satisfies a {\bf WLSI} with  function $\beta_{WL}$, then it verifies
 a classical Poincar\'e inequality if and only if there exists
$c_1,c_2>0$ such that for $s$ small enough,
$$
\beta_{WL}(s) \leq c_1 \log(c_2/s).
$$
\end{eprop}

\begin{eproof}
For the first statement, first note that $\beta_{WP}$ is non-increasing. Then
Theorem~\ref{thm-cap} and Remark~\ref{rem-eq} imply that
$$
\frac{\frac{{\mu(A)}}{2}\log\PAR{1+\frac{1}{2\mu(A)}}}{\beta_{WL}\PAR{\frac{\mu(A)}{2}
\log\PAR{1+\frac{1}{2\mu(A)}}}}\leq Cap_\mu(A).
$$
This means that
$$
\frac{12\mu(A)}{\beta_{WP}(\mu(A))}\leq Cap_\mu(A),
$$
where $\beta_{WP}$ is defined by \eqref{eq-dd}, the result holds using Theorem~2.2 of
\cite{ca-ba-ro2}.

To prove the second statement we use the same argument (replacing Theorem \ref{thm-cap} by Theorem
\ref{thm-wls}) and the fact that there exist constants $A,A',s_0>0$ such that
\begin{equation}
\label{eq-A}
\forall s\in(0,s_0),\quad A'\frac{s}{\log\PAR{1/s}}\leq \phi^{-1}(s)\leq A \frac{s}{\log\PAR{1/s}},
\end{equation}
where $\phi(s)=s \, \log(1+e^2/s)$. Then $\mu$
satisfies a {\bf WLSI} with function $\beta_{WL}$ defined by \eqref{eq-ddd}. Note that $\beta_{WL}$  is non-increasing.

Finally, the last two results prove that $\beta_{WL}(s) \leq c_1 \log(c_2/s)$ for $s$ enough is equivalent to classical
Poincar\'e inequality.
\end{eproof}

\begin{erem}
\begin{itemize}
\item It is interesting to remark that when considering the usual derivation ``Logarithmic Sobolev
inequality implies Poincar\'e inequality'' by means of test function $1+\epsilon g$ and
$\epsilon\to0$, we get a worse result: a weak logarithmic Sobolev inequality  with function $\beta$
implies a weak Poincar\'e inequality  with the same function $\beta$, whereas the result
of the proposition~\ref{prop-dd} gives a better result.
\item As a byproduct, we get that any Boltzman's measure
(with a locally bounded potential) satisfies some {\bf WLSI} if Ricci($M$) is bounded from
below (see \cite{r-w}).
\item Finally the above proof shows that we obtain the best function (up to multiplicative
constants) for {\bf WPI} or {\bf WLSI} as soon as we have the best function for the other. In
particular we recover the good functions for the examples~\ref{ex1}.
\end{itemize}
\end{erem}
\medskip

\subsection{Link with super Poincar\'e inequalities}

Let us recall the definition of the super Poincar\'e inequality introduced by Wang in \cite{w1}.

\begin{edefi}
\label{def-sp} We say that the measure $\mu$ satisfies  a super Poincar\'e inequality, {\bf SPI}, if
there exists a non-increasing function $\beta_{SP} : [1,+\infty)\rightarrow\dR^+$, such that for
all $s \geq 1$ and any function $f\in H^1(M,\mu)$,
\begin{equation}
\tag{{\bf SPI}} \int f^2d\mu\leq\beta_{SP}(s)\int\ABS{\nabla f}^2d\mu+s \, \left(\int \ABS{f}
d\mu\right)^2 .
\end{equation}
 \end{edefi}

Note that as for {\bf WLSI} in Remark~\ref{rema},  for the {\bf SPI} what is important is
the behaviour of $\beta$ near~$\infty$. As for Proposition~\ref{prop-dd} we can now relate {\bf WLSI} and {\bf SPI}.

\begin{eprop}
\label{prop-super}
Suppose that $\mu$ satisfies a {\bf WLSI} with function $\beta_{WL}$. Assume that $\beta_{WL}$ verifies
that
$
x\mapsto {\beta_{WL}\PAR{\frac{\log(x/2)}{2x}}}/{\log(x/2)}
$
is non-increasing on $(2,\infty)$.

Then $\mu$ satisfies a {\bf SPI} with function $\beta_{SP}$ given by
\begin{equation}\label{eq-sp1}
\forall t\geq 2e,\qquad\beta_{SP}(t)= 2\frac{\beta_{WL}\PAR{\frac{\log(t/2)}{2t}}}{\log(t/2)},
\end{equation}
and constant on $[1,2e)$.
\end{eprop}

\begin{eproof}
If $\mu$ satisfes a {\bf WLSI} then one obtains by  Theorem~\ref{thm-cap} and Remark~\ref{rem-eq}:
\begin{equation}
\label{eq-prof}
\frac{\frac{\mu(A)}{2}\log\PAR{1+\frac{1}{2\mu(A)}}}{\beta_{WL}\PAR{\frac{\mu(A)}{2}
\log\PAR{1+\frac{1}{2\mu(A)}}}}\leq Cap_\mu(A),
\end{equation}
for any $A\subset M$, with $\mu(A)\leq1/2$. Finally the function
$
t \mapsto t\, {\beta_{WL}\PAR{\frac{\log(t/2)}{2t}}}/{\log(t/2)}
$ is clearly non decreasing for $t\geq 2e$, then Corollary~6 of~\cite{ca-ba-ro3} gives the result.
\end{eproof}

Note that the last proposition is not entirely satisfying, we hope that {\bf WLSI} is equivalent to {\bf SPI} via
a measure-capacity measure criterion.


\subsection{Link with general Beckner  inequalities}

\begin{edefi}
\label{def-bec}
Let $T : [0,1]\rightarrow\dR^+$, be a non-decreasing function, satisfying in addition
$x\mapsto T(x)/x$ is non-increasing on $(0,1]$.

We say that a measure $\mu$ satisfies  a general Beckner inequality, {\bf GBI}, with function $T$
 if for all function $f\in H^1(M,\mu)$,
\begin{equation}
\tag{{\bf GBI}} \sup_{p\in(1,2)} \, \frac{\int f^2 d\mu \, - \, \left(\int |f|^p d\mu\right)^{\frac
2p}}{T(2-p)} \, \le\, \int\ABS{\nabla f}^2d\mu.
\end{equation}
\end{edefi}

Note that our hypotheses imply that
$$
\forall x\in[0,1],\qquad T(1)x \leq T(x) \leq T(1).
$$
 The two extremal cases correspond
respectively to the Poincar\'e inequality ($T$ is constant, $T(x)=T(1)$) and the logarithmic Sobolev inequality
($T(x)=T(1)x$). The intermediate cases $T(x)=x^a$ for $0\leq a \leq 1$, have been introduced and studied
in \cite{LO00}, while a study of general $T$ is partly done in \cite{ca-ba-ro}. Also note that (up
to multiplicative constants) the interesting part of $T$ is its behaviour near 0, that is we can
always define $T$ near the origin and then take it equal to a large enough constant. Recall finally that the usual Beckner inequality concerns $T(x)=x$ and was introduced by Beckner  to get quantitative information on an interpolation between Poincar\'e's inequality and logarithmic Sobolev inequality for the Gaussian measure, see \cite{1989beckner}.

In \cite{ca-ba-ro} Theorem 10 and Lemma 9, it is shown that (up to a multiplicative
constant 3) {\bf GBI} is equivalent to a measure-capacity. More precisely, the inequality measure-capacity~\eqref{def-mc}
with the function
\begin{equation}
\label{eq-gb}
\gamma(u)=T\left(\frac{1}{\log\left(1+\frac 1u\right)}\right),
\end{equation}
for $u$ small enough implies a  {\bf GBI} with the function $20T$. And  {\bf GBI} implies a measure-capacity inequality with the function
$6\gamma$ defined on~\eqref{eq-gb}.   We thus obtain:

\begin{eprop}\label{prop-becwl}
Assume that $\mu$ satisfies a {\bf WLSI} with function $\beta_{WL}$. Let
\begin{equation}
\label{eq-wlbec}
\forall t\in(0,1],\qquad T(t)= t\beta_{WL}\PAR{\frac{1}{4te^{1/t}}},
\end{equation}
and assume that $T$ non-decreasing $(0,t_a]$ for some $t_a\in(0,1]$.
Then the measure $\mu$ satisfies a {\bf GBI} with function $20T$.

Conversely assume  that $\mu$ satisfies a {\bf GBI} with function $T$ and constant $c$, then $\mu$
satisfies a {\bf WLSI} with function $\beta_{WL}$ given by
\begin{equation}\label{eq-becwl}
\beta_{WL}(s)=C\, T\left(C'\frac{1}{\log(1/s)}\right) \, \log(1/s),
\end{equation}
for $s$ small enough and some  constants $C,C'$.
\end{eprop}

\begin{eproof}
Assume that $\mu$ satisfies a {\bf WLSI} with function $\beta_{WL}$.
Using Theorem~\ref{thm-cap} and Remark~\ref{rem-eq} one has inequality~\eqref{eq-prof}.
Using the fact that
$$
\forall x\in(0,1],\qquad \log\PAR{1+\frac{1}{2x}}\geq\frac{1}{2}\log\PAR{1+\frac{1}{x}},
$$
one obtains that inequality~\eqref{eq-prof} implies that the function $T$ defined on~\eqref{eq-wlbec} satisfies
a measure-capacity inequality. The function $x\mapsto T(x)/x$ is non-increasing and
due to  the fact that $T$ is non-decreasing by hypothesis, then
  Theorem~10 and Lemma~9 of  \cite{ca-ba-ro} prove that $\mu$ satisfies a
{\bf GBI} of function $T$.

To prove the second statement we need also Theorem~10 and Lemma~9 of \cite{ca-ba-ro}, Theorem~\ref{thm-wls} and
inequality~\eqref{eq-A}.
\end{eproof}

\begin{eex}
Note that if the function $T$ defined on~\eqref{eq-wlbec} is non-decreasing near 0 then one can prove
that $\beta_{WL}(s)\leq c_1\log(c_2/s)$ for $s$ small enough and some constants $c_1,c_2>0$. Then by
Proposition~\ref{prop-dd},
$\mu$ satisfies a  Poincar\'e inequality. The last proposition can be applied only for measures satisfying
a Poincar\'e inequality.
\end{eex}

\subsection{Link with an other weak logarithmic Sobolev inequality}
The next inequality is useful to control the decay in entropy of the semigroup. It will by used in
Theorem~\ref{thm-resls}.

\begin{ethm}\label{thm-wlsvar}
If $\mu$ satisfies a {\bf WLSI} with function $\beta_{WL}$, then $\mu$ satisfies for any function $f\in H^1(M,\mu)$ and any $u>0$ small enough,
\begin{equation}
\ent{\mu}{f^2} \leq \beta_{SWL}(u)\int\ABS{\nabla f}^2d\mu+\sqrt{3}u \, \left(\varf{\mu}(f^2)\right)^{\frac
12} \, ,
\end{equation}
with $$\beta_{SWL}(u)= 16 \beta_{WL}\left(\frac{\kappa \, u^3}{\log^6(1/u)}\right)
$$ for some universal constant $\kappa$ and $u$ small enough.
\end{ethm}

\begin{eproof}
According to Theorem~\ref{thm-cap} and Remark~\ref{rem-eq} we know that for every $A\subset M$ such that
$\mu(A)\leq1/2$,
$$
Cap_\mu(A) \geq
\frac{\frac{\mu(A)}{2}\log\PAR{1+\frac{1}{2\mu(A)}}}{\beta_{WL}\PAR{\frac{\mu(A)}{2}\log\PAR{1+\frac{1}{2\mu(A)}}}}
\geq
\frac{\frac{\mu(A)}{2k}\log\PAR{1+\frac{e^2}{\mu(A)}}}{\beta_{WL}\PAR{{\mu(A)\log\PAR{1+\frac{e^2}{\mu(A)}}}}}
$$
for $k=\log(1+2e^2)/\log(2)$ using $k \log(1+y/2)\geq \log(1+e^2 y)$ for $y\geq 2$ and that
$\beta_{WL}$ is non-increasing. Hence we are in the situation of Theorem \ref{thm-wls} with $\beta=
2k \, \beta_{WL}$.

Note that we may assume that $f$ is non-negative.
\smallskip


Let $\Omega_0\subset M$, it will be fixed latter. Indeed the first quantity we have to control is $\int_{\Omega_0} \, F_+^2 h d\mu$ which is less
than
$$
\left(\int_{\Omega_0} h^2 d\mu\right)^{\frac 12} \, \left(\int F_+^4 d\mu\right)^{\frac 12}\, .
$$
We thus have to bound
\begin{eqnarray*}
X_0 & := & \sup\{\int_{\Omega_0} h^2 d\mu; \, h\geq 0, \, \int e^h d\mu \leq 1+e^2\} \\ & = &
\sup\{\int_{\Omega_0} h^2 d\mu; \, h\geq 0, \, \int_{\Omega_0} e^h d\mu \leq e^2 + \mu(\Omega_0)\}
\, ,
\end{eqnarray*}
(see \cite{barthe-roberto} Lemma 6 for the latter equality). But $\phi(x)=(1+\log^2(x)) \,
\1_{x\geq e} \, + \, \frac 2e x  \, \1_{x<e}$ is concave and non-decreasing on $\dR_+$. It follows that
\begin{eqnarray*}
\phi\left(\frac{e^2+\mu(\Omega_0)}{\mu(\Omega_0)}\right) & \geq & \int_{\Omega_0} \phi(e^h) \,
\frac{d\mu}{\mu(\Omega_0)} \\ & \geq & \int_{\Omega_0} \left((1+h^2) \, \1_{h\geq 1}\right) \,
\frac{d\mu}{\mu(\Omega_0)}\\ & \geq & \int_{\Omega_0} h^2 \, \frac{d\mu}{\mu(\Omega_0)} \, - \,
\int_{\Omega_0} h^2 \, \1_{h< 1} \, \frac{d\mu}{\mu(\Omega_0)}\\& \geq & \int_{\Omega_0} h^2 \,
\frac{d\mu}{\mu(\Omega_0)} \, - \, 1 \, ,
\end{eqnarray*}
so that $$X_0 \, \leq \, \mu(\Omega_0) \, \left(2 +
\log^2\left(1+\frac{e^2}{\mu(\Omega_0)}\right)\right) \, := \, \psi(\mu(\Omega_0)) \, .$$
Once again onlys small values of $s$ are challenging, consider then $s\leq 1$.  We can mimic now the proof of Theorem \ref{thm-wls} , briefly  we define $c$ by
$$
c=\inf\BRA{t\geq 0,\quad \psi( \mu(F_+^2>c))\leq s},
$$
and then we choose $\Omega_0$ such that
$\BRA{F^2_+>c}\subset  \Omega_0\subset \BRA{F^2_+\geq c}$ and
$
\psi( \mu(\Omega_0))=s^a,
$
for some $a>0$.  This
choice being possible since $\psi$ is increasing on $[0,1/2[$, the maximal possible $s$ being
greater than 1.

Then and obtain
\begin{multline}
\label{eq-tata}
\disp\sup\BRA{\int F_+^2hd\mu;\,\,h\geq0,\,\,\int e^hd\mu\leq e^2+1} \,  \leq   \,  \disp
\sqrt{s^a} \, \left(\int F_+^4 d\mu\right)^{\frac 12} \  +  s \frac{c}{1-\rho}\\
 + \frac{\beta_{WL}(s)}{\rho(1-\sqrt{\rho})^2} \int\ABS{\nabla
F_+}^2d\mu \, .
\end{multline}
It remains to estimate $c$. Note that there exists an universal constant $\theta$ such that
$\psi^{-1}(x)\geq \theta \, x/\log^2(1+\frac{e^2}{x})$ . It follows using this two inequalities
$$
\theta \,\frac{s^a}{\log^2(1+\frac{e^2}{s^a})} \leq \mu(\Omega_0)
$$
and by Markov inequality
$$
 \mu(\Omega_0) \leq \frac{\int F_+^2 d\mu}{c} \leq
\frac{\left(\int F_+^4 d\mu\right)^{\frac 12}}{c} \, ,
$$
so that choosing $a=2/3$ and $\rho=1/4$ we
finally obtain
\begin{multline}
\label{eq-tatata}
\disp\sup\BRA{\int F_+^2hd\mu;\,\,h\geq0,\,\,\int e^hd\mu\leq e^2+1} \leq  {s^{\frac 13}} \,
\left(1+\frac{4}{3\theta} \log^2(1+\frac{e^2}{s^{2/3}})\right) \left(\int F_+^4
d\mu\right)^{\frac 12} \\  +  16 \beta_{WL}(s) \int\ABS{\nabla F_+}^2d\mu \, .
\end{multline}
The same inequality for $F_-$ and the elementary $\sqrt{a}+\sqrt{b}\leq \sqrt{2} \, \sqrt{a+b}$
yield, since there exists an universal constant $\theta'$ such that the inverse function of $s
\mapsto \sqrt{2} s^{\frac 13} \left(1+\frac{4}{3\theta} \log^2(1+\frac{e^2}{s^{2/3}})\right)$ is
greater than $u \mapsto \theta' \, u^3 / (\log^6(1/u))$ for $u>0$ small enough,
\begin{equation}\label{eq-enfin}
\ent{\mu}{f^2} \leq 16 \beta_{WL}\left(\frac{\theta' \, u^3}{1+\log^6(1+\frac{e^2}{u^2})}\right)
\int\ABS{\nabla f}^2d\mu \, + \, u \, \left(\int (f-m)^4 d\mu\right)^{\frac 12} \, .
\end{equation}
Since we have assumed that $f$ is non-negative, a median of $f^2$ is $m^2$, and $(f-m)^4 \leq
(f^2-m^2)^2$. Finally, if $M$ denotes the mean of $f^2$, $$\int \left((f^2-M)-(m^2-M)\right)^2 d\mu
=\varf{\mu}(f^2)+(m^2-M)^2$$ and since $m^2-M$ is a median of $f^2-M$, provided $m^2-M\geq 0$
$$\varf{\mu}(f^2)\geq \int (f^2-M)^2 \, \1_{f^2-M\geq m^2-M} \, d\mu \geq \frac 12 \, (m^2-M)^2$$
while if $m^2-M\leq 0$ $$\varf{\mu}(f^2)\geq \int (f^2-M)^2 \, \1_{f^2-M\leq m^2-M} \, d\mu \geq
\frac 12 \, (m^2-M)^2 \, .$$ We thus finally obtain  $$\int (f-m)^4 d\mu \leq 3 \,
\varf{\mu}(f^2)$$ and the proof is completed.
\end{eproof}

One may of course derive other weak logarithmic Sobolev inequalities by this method, such
inequalities as well as further applications will be treated elsewhere. We will apply this theorem on the
section~\ref{section-5} for the decay to the equilibrium of the semigroup.

\section{Convergence of the associated semigroup}
\label{section-5}

In this section we shall study entropic convergence for the associated semi-group. Namely we
assume that $(\Pt)_{t\geq 0}$ is a ``nice'' diffusion $\mu$ symmetric semi-group. Here by ``nice''
we mean that $(\Pt)_{t\geq 0}$ is the semi-group associated to a non-explosive diffusion process
on some Polish space admitting a ``carr\'e du champ''. For a precise framework we refer to
\cite{cat4} Section 1.1. Roughly speaking, these assumptions allow us to give a rigorous meaning
to all computations below.
\smallskip

Let $h$ be a bounded density of probability with respect to the measure $\mu$. The two results of
this section connect the decay of the entropy with the infinite norm of $h$. More precisely, using
the {\bf WLSI} we will compute the function $C(t,\NRM{h}_\infty)$ such that for all $t>0$,
 $$
\ent{\mu}{\Pt h} \leq C(t,\NRM{h}_\infty).
 $$
 We will give here conditions under which $C(t,\NRM{h}_\infty)\rightarrow 0$ when $t$ goes to $\infty$.

\medskip

The first result connects the decay of the entropy with the oscillation of $h$, one gets:
\begin{eprop}\label{prop-decay1}
Let  $\mu$ satisfies a {\bf WLSI} with function $\beta_{WL}$ and let $h\geq 0$, bounded with $\int
hd\mu=1$. Then  for any $\varepsilon>0$ and for $t$ large enough we get:
\begin{equation}
\label{eq-d1} \ent{\mu}{\Pt h}\le(e^{-1}+\varepsilon)\,\xi_{\varepsilon}(t) \, {\bf Osc}^2(\sqrt
h)
\end{equation}
where $\xi_\varepsilon(t)$  is given by
$$
\xi_{\varepsilon}^{-1}(r)=-{1\over2}\beta_{WL}(r)\log\PAR{\frac{ r}{\varepsilon}},
$$ for $r$ small enough.

Conversely, if there exists a decreasing function $\xi$ such that, for any bounded $h\geq 0$, with
$\int h d\mu=1$  we have
$$
\forall t>0,\qquad\ent{\mu}{\Pt h}\le\xi(t) \, {\bf Osc}^2(\sqrt h),
$$
then $\mu$ satisfies a {\bf WLSI} with function $\beta_{WL}(t)=\psi^{-1}(t)$ where $\psi(t)= 2
\sqrt{ 2 \, \xi(t)}$. In particular if $\xi(t)\leq c e^{-\alpha t}$, for some $\alpha>0$, the
measure $\mu$ satisfies a Poincar\'e inequality.
\end{eprop}
\begin{eproof}
We start with the direct part. Denote $I(t)=\ent{\mu}{\Pt h}$. Then $I'(t)= - \, \frac 12 \, \int
\frac{|\nabla \Pt h|^2}{\Pt h} d\mu$, thus the weak logarithmic Sobolev inequality yields
\begin{eqnarray*}
I'(t)&\le & -{2\over \beta_{WL}(r)}I(t)+{2r\over\beta_{WL}(r)}{\bf Osc}^2(\sqrt{\Pt h}).
\end{eqnarray*}
Using Gronwall's lemma yields
$$\ent{\mu}{\Pt h}\le \inf_{r>0}\left\{ r\sup_{s\in[0,t]}{\bf Osc}^2(\sqrt{\Ps h}) + e^{-2t/\beta_{WL}(r)}\ent{\mu}{h}\right\}.$$
We may now use ${\bf Osc}^2(\sqrt{\Pt h})\leq {\bf Osc}^2(\sqrt{h})$ and $\ent{\mu}{h} \leq 1/e \,
{\bf Osc}^2(\sqrt{h})$ as we quoted in Remark \ref{rema} and finally choose $r$ such that
$r=\varepsilon \, e^{-2t/\beta_{WL}(r)}$ (which is optimal up to constants) to get the result.

Let us prove the second statement.
Denote $f=\sqrt{h}$. According to \cite{cat4} (2.5) with $\alpha_1=-1$ and $\alpha_2=2$ it holds
\begin{equation}\label{eq-lsrobust}
\ent{\mu}{h} \le t \, \int\ABS{\nabla f}^2 d\mu + 2 \, \log \, \int f \, \Pt h \, d\mu \, .
\end{equation}
But
\begin{eqnarray*}
\int f \, \Pt h \, d\mu & = & \int f \, \left(1+(\Pt h - 1)\right) \, d\mu \\ & \leq & 1 \, + \,
\int (f - \int f d\mu) \, (\Pt h - 1) \, d\mu \\ & \leq & 1 \, + \, {\bf Osc}(f) \, \int |\Pt h -
1| \, d\mu \\ & \leq & 1 \, + \, {\bf Osc}(f) \, \sqrt{2 \, \ent{\mu}{\Pt h}} \\ & \leq & 1 \, + \,
\sqrt{2 \xi(t)} \, {\bf Osc}^2(f) \, ,
\end{eqnarray*}
where we used successively $\int f d\mu \leq 1$ , Pinsker inequality and the hypothesis. It remains
to use $\log(1+a) \leq a$ to get the first result.
The particular case follows from Proposition \ref{prop-dd}.
\end{eproof}

The previous result is the exact analogue of Theorem 2.1 in \cite{r-w} for {\bf WPI}. The converse
statement (Theorem 2.3 in \cite{r-w}) is remarkable in the following sense: it implies in
particular that any exponential decay ($\varf{\mu}(\Pt f) \leq c e^{-\alpha t} \Psi(f-\int f
d\mu)$) for any $\Psi$ such that $\Psi(af)=a^2 \Psi(f)$ (in particular $\Psi(f)={\bf Osc}^2(f)$)
implies a (true) Poincar\'e inequality. This result is of course very much stronger than the usual
one involving a $\dL^2$ bound. Its proof lies on the fact that $t \mapsto \log (\int (\Pt f)^2
d\mu)$ is convex. This convexity property (even without the log) fails in general for the relative
entropy (Bakry-Emery renowned criterion was introduced for ensuring such a property).
Actually a similar statement for the entropy is false.

\smallskip

Not that the previous result is only partly satisfactory for the convergence of the entropy.
Indeed recall that for a density of probability $h$, the following holds
$$\varf{\mu}(\sqrt{h})\leq \ent{\mu}{h}\leq \varf{\mu}(h)$$ so that a weak Poincar\'e inequality
implies for $t>0$
$$
\ent{\mu}{\Pt h} \leq \xi^{WP}_{\varepsilon} (t)\, (1+\varepsilon){\bf Osc}^2(h) \, ,
$$
whereas our WLSI implies
$$
\ent{\mu}{\Pt h} \leq \xi^{WLS}_{\varepsilon}(t) \, (e^{-1}+\varepsilon){\bf Osc}^2(\sqrt{h}),
$$
so that for small time, the WLSI  furnishes better bounds than a weak Poincar\'e inequality (and
justifies the use of LSI for this kind of evaluation), though the rate of convergence is not the
expected one.


\medskip

In order to correct this unsatisfactory point, at least when a Poincar\'e inequality holds, and
always for bounded density $h$, we will make use of the other weak logarithmic Sobolev inequality
stated in Theorem \ref{thm-wlsvar}. Indeed, another way to control entropy decay was introduced in
\cite[Theorem 1.13]{gui-ca}. It was proved there that a Poincar\'e inequality (with constant
$C_P$) is equivalent to a \underline{restricted} logarithmic Sobolev inequality
$$
\ent{\mu}{h} \leq C \, (1+\log(\parallel h\parallel_{\infty})) \, \int\frac{|\nabla h|^2}{h} \, d\mu
$$
 for all bounded density of probability $h$, where the constant $C$ only depends
on $C_P$. It follows that
$$\ent{\mu}{\Pt h} \leq e^{- \, \frac{t}{C(1+\log(\parallel h\parallel_{\infty}))}} \,
\ent{\mu}{h}$$ for such an $h$.
\smallskip


We shall describe below one result in this direction for {\bf WLSI}, using
Theorem~\ref{thm-wlsvar} and Poincar\'e inequality.

\begin{eprop}\label{thm-resls}
Let $\mu$ be a probability measure satisfying a {\bf WLSI} with function $\beta_{WL}$ and the usual
Poincar\'e inequality with constant $C_P$. Let $\beta_{SWL}$ be the function defined in
Theorem~\ref{thm-wlsvar}. Then for all $f\in H^1(E,\mu)$,
$$
\ent{\mu}{f^2} \leq A(C_P,\parallel
f\parallel_{\infty}) \, \int\ABS{\nabla f}^2d\mu
$$
 where
$$
A(C_P,\NRM{f}_{\infty})=\inf_{u\in(0,s_0]}\BRA{\beta_{SWL}(u)+u\sqrt{3C_P}\NRM{f}_\infty^2}.
$$
Here $s_0$ is any positive number such that $\beta_{SWL} $ is defined by the formula in Theorem
\ref{thm-wlsvar} for $s \leq s_0$ and then extended by $\beta_{SWL}(s)=\beta_{SWL}(s_0)$ for
$s\geq s_0$.

As a consequence, for all $t\geq0$,
$$
\ent{\mu}{\Pt h} \leq e^{- \, t \, / \, A\PAR{C_P,\sqrt{\parallel h\parallel_{\infty}}}} \, \ent{\mu}{h}
$$
 for any bounded density of probability $h$.
\end{eprop}

\begin{eproof}
Due to homogeneity we may assume that $\int\ABS{\nabla f}^2d\mu=1$ (if it is 0 the result is
obvious). But since $\mu$ satisfies a Poincar\'e inequality
$$
\varf{\mu}(f^2) \leq 4 C_P \, \int
f^2 \, |\nabla f|^2 d\mu \leq 4 C_P \parallel f^2\parallel_{\infty} \, ,
$$
 so that
$
\ent{\mu}{\Pt f^2}
\leq \beta_{SWL}(u) + 2u \parallel f\parallel_{\infty}^2 \, \sqrt{3C_P} \, .
$
\end{eproof}

Note now that the previous entropic decay is always better for small time. Indeed if
$$t\le {C_P A(C_P,\parallel h\parallel_{\infty}^{\frac 12})\over A(C_P,\parallel h\parallel_{\infty}^{\frac 12})-C_P}\log\left({\varf{\mu}(h)\over\ent{\mu}{ h}}\right)$$
then the entropic decay obtained by Proposition \ref{thm-resls} is better than the estimate with
Poincar\'e inequality.

\begin{eex}\label{cor-rlslo}
Let $\alpha \in [1,2]$ and $d\mu_\al(t)=Z_\al e^{-|t|^{\al}}dt,\; t\in \dR$ where $Z_\al$ is a normalization constant. Using Example~\ref{ex1} and Proposition~\ref{prop-becwl} one obtains  that $\mu_\al$ satisfies a {\bf GBI} with $T(x)=C \, x^{\frac{2\alpha - 2}{\alpha}}$ for $x\in(0,1)$. Then one can find $C(\alpha),C'(\alpha)>0$  such that for all bounded density of probability $f$,
 $$
 \ent{\mu}{f^2} \leq C(\alpha) \,
\left(1+\log^{(2/\alpha)-1}(\parallel f\parallel_{\infty})\right) \, \int\ABS{\nabla f}^2d\mu \,.
$$
As a consequence, for all $t\geq0$,
$$
\ent{\mu}{\Pt h} \leq e^{- \, t \, / \, C'(\alpha) \,
\left(1+\log^{(2/\alpha)-1}(\parallel h\parallel_{\infty})\right)} \, \ent{\mu}{h},
$$
 for any bounded density of probability $h$.
\end{eex}

It seems very unlikely that one can derive such a result from a direct use of Proposition~\ref{prop-decay1}.
As noticed in \cite{gui-ca}, these restricted logarithmic Sobolev inequalities (restricted to the ($\Pt$
stable) $\dL^{\infty}$ balls) can be used to obtain modified (or restricted) transportation
inequalities. We recall below a result taken from section 4.2 in \cite{gui-ca}. If $\nu=h \mu$ is a
probability measure, it can be shown
\begin{equation}\label{eq-trans}
W_2^2(\nu,\mu) \, \leq \,  \eta(0) \, \ent{\mu}{h} \, + \, \int_0^{+\infty} \, \eta'(t) \,
\ent{\mu}{\Pt h} \, dt \, ,
\end{equation}
where $\eta$ is a non-decreasing positive function such that $\int (1/\eta(t)) dt=1$, and $W_2$
denotes the (quadratic) Wasserstein distance between $\nu$ and $\mu$. We may take here
$$\eta(t)=2 \, A(C_P,\parallel h\parallel_{\infty}^{\frac 12}) \,
e^{\frac 12 \,  t / A(C_P,\parallel h\parallel_{\infty}^{\frac 12})}$$ which yields
\begin{equation}\label{eq-trans2}
W_2(\nu,\mu) \, \leq \, D \, (1+A^{\frac 12}(C_P,\parallel h\parallel_{\infty}^{\frac 12})) \,
\sqrt{\ent{\mu}{h}} \, .
\end{equation}
In the Latala-Oleszkiewicz situation, we recover, up to the constants, Theorem 1.11 in
\cite{gui-ca}.

Using Marton's trick, \eqref{eq-trans2} allows us to obtain a concentration result (a little bit
less explicit than the one obtained via {\bf GBI} in Proposition 29 of \cite{ca-ba-ro}) namely
there exist $r_0$ and $\sigma$ such that if $\mu(A)\geq 1/2$ and $A_r^c=\{x,d(x,A)\geq r\}$ one has
$$r-r_0 \, \leq \, \sigma \, A^{\frac 12}(C_P,(1/\mu^{1/2}(A_r^c))) \, \sqrt{\log(1/\mu(A_r^c))} \, .$$
In the Latala-Oleszkiewicz situation, we recover up to the constants, the same concentration
function as $\mu_\al$, showing that our restricted logarithmic Sobolev inequality is (up to the constants)
optimal. Note that another way to get the concentration result is to use the modified logarithmic Sobolev
(and transportation) inequalities discussed in \cite{GGM1,GGM2}.\par\vspace{5pt}

Let us finally note that even if the results obtained by the WLSI are always efficient in the
regime between Poincar\'e and Gross inequality, it relies on the crucial assumption that $h$ is a
bounded density. The goal of the next section is to get rid of this assumption.

\section{Convergence to equilibrium for diffusion processes}
\label{sec-dif}

In this section we shall discuss the rate of convergence to equilibrium for particular diffusion
process, both in total variation and in entropy. The main difference between the previous section
is that we do not  assume that the initial law of the diffusion processes has a density of
probability with respect to symmetric measure $\mu$. The initial entropy is not necessarily
finite.

For simplicity we only consider the case when $M=\dR^n$ and $\mu=e^{-2V}dx$. Hence our diffusion
process is given by the stochastic differential equation
\begin{equation}
\label{eds}
dX_t \, = \, dB_t \, -  \, (\nabla V)(X_t)dt \quad , \quad Law(X_0)=\nu
\end{equation}
where $B_.$ is a standard Brownian motion. We assume that $V$ is $C^3$ and that there exists some
$\psi$ such that $\psi(x) \rightarrow +\infty$ as $|x| \rightarrow +\infty$ and $\frac 12 \,
\Delta \, \psi \, - \, \nabla V.\nabla \psi$ is bounded from above. This assumption ensures the
existence of an unique non explosive strong solution for \eqref{eds}. If $\nu=\delta_x$ we will
denote by $X_t^x$ the associated process (cf e.g. \cite{Ro99}).

A remarkable consequence of Girsanov theory (see \cite{Ro99} in our situation) is that with our
assumptions, for all $\nu$ and all $t>0$ the law of $X_t$ denoted by $\Pt \nu$ is absolutely
continuous with respect to $\mu$, its density will be denoted by $h_t$. Of course if $\nu=h \mu$,
$\Pt \nu= (\Pt h) \mu$ and $\mu$ is a reversible measure.

In particular $\Pt \nu = (\P_{t-u} h_u) \mu$, and the rate of convergence of $\Pt \nu$ towards
$\mu$ can be studied by using the semigroup properties only. In the sequel we shall make the abuse
of notation $\Pt \nu = (\P_{t-u} h_u)$ i.e. we shall abusively identify the measure with its
density. What we need to understand is thus the behavior of $\Pt h$, where $h$ is a density of
probability (in a sense it is $\Pt f^2$ rather than $\Pt f$ which is interesting).
\smallskip

Of particular interest is the case when
\begin{equation}\label{condition}
|\nabla V|^2(x) \, - \, \Delta V(x) \, \geq \, - C_{min} \, > - \infty
\end{equation}
 for a nonnegative
$C_{min}$ since in this case one can show (see \cite[Theorem 3.2.7]{Ro99}) that $\ent{\mu}{\Pt
\delta_x}$ is finite for all $t>0$. Actually the proof of Royer can be used in order to get the
following more general and precise result
\begin{eprop}\label{prop-royer}
With the previous hypotheses
\begin{multline}\label{eq-royer1}
\int \Pt \delta_x \, \log_+^p(\Pt \delta_x) \, d\mu \, \leq \\ \, 4^{p-1}
\left(V_+^p(x)+\left(\frac{C_{min} t}{2}\right)^p+\left(\frac n2 \, \log(\frac{1}{2\pi
t})\right)^p+e^{V(x)+p(\log p -1)+\frac 12 C_{min}t}\right)
\end{multline}
for all $t\in ]0,1/2\pi[$ and $p\geq 1$.

If in addition
\begin{equation}\label{eq-initial}
V_+(y) \leq D (V_+(x) + |y-x|^2 + D')
\end{equation}
 for some $D>0$, $D'$ and all pair $(x,y)$, then for all $t\in ]0,1/2D\wedge 1/2\pi[$
\begin{equation}\label{eq-royer2}
\int \Pt \delta_x \, \log_+^p(\Pt \delta_x) \, d\mu \, \leq \, 4^{p-1} \left((1+D^p) \,
(V_+(x)+D')^p+\left(\frac{C_{min} t}{2}\right)^p+\left(\frac n2 \, \log(\frac{1}{2\pi
t})\right)^p\right) \, .
\end{equation}

In particular, if $\int e^{V_+} d\nu := M <+\infty$,
\begin{equation}
\label{eq-nuc}
\left(\int \Pt \nu \, \log_+^p(\Pt \nu) \,
d\mu\right)^\frac 1p \, \leq \, p \, C(\nu,t_0)
\end{equation}
 for all $t\geq t_0 >0$, where $C(\nu,t_0)$ only
depends on $t_0$, $M$, $C_{min}$ ($\lambda$) and the dimension. If in addition~\eqref{eq-initial}
holds, it is enough to assume that $\int e^{\lambda \, V_+} d\nu := M <+\infty$ for some $\lambda
0$.
\end{eprop}
\begin{eproof}
Let $$F = \exp \, \left(V(x)-V(W_t)-\frac 12 \, \int_0^t \left(|\nabla V|^2-\Delta V\right)(W_s)
ds\right) \, ,$$ where $W$ is a Brownian motion starting from $x$. Recall that $F$ is a density of
probability (with our hypotheses). If $I(t)=\int \Pt \delta_x \, \log_+^p(\Pt \delta_x) \, d\mu$ we
may use the argument in \cite[Theorem 3.2.7]{Ro99} and the convexity of $u \mapsto u^p$ in order to
get
$$I(t) \leq \dE \left(F \, 4^{p-1} \left(V_+^p(x)+(V(W_t)-\frac{1}{2t} |W_t-x|^2)_+^p +
(C_{min}t/2)^p + \frac n2 \log(\frac{1}{2\pi t})\right)^p \right) \, .$$ The first statement
follows easily bounding $(V(W_t)-\frac{1}{2t} |W_t-x|^2)_+$ by $D (V(W_t)+D')_+$ and $u^p e^{-u}$
by $p^p \, e^{-p}$. The second one is immediate since \eqref{eq-initial} allows us to bound the
same term by $V_+(x)$ for $t$ small enough.

The last statements are obtained by using two arguments. First $u^p \leq p! \, e^u$ (or $u^p \leq
p! (1/\lambda)^p e^{\lambda p}$), so that for a given $t$ the result follows from $(p!)^\frac 1p
\leq c p$. The second one is standard, namely $t \mapsto \int \Pt h \, \log_+^p \Pt h \, d\mu$ is
non-increasing.
\end{eproof}
\smallskip

We shall come back to the condition \eqref{eq-initial} later on. Note however that such a
condition is trivially verified for $V(x)=|x|^\gamma$, $0<\gamma\le 2$.

\subsection{Rate of convergence for the relative entropy}

\begin{ethm}\label{thm-entdecay}
Let $d\mu=e^{-2V} dx$ be a probability measure which satisfies a {\bf WLSI} with function
$\beta_{WL}$ and let $\xi$ be defined as in~\eqref{eq-d1} of Proposition~\ref{prop-decay1}. Assume that
\eqref{condition} holds and let $\nu$ be a probability measure such that~\eqref{eq-nuc} holds.

Then for all $1\geq \varepsilon >0$ and all $k>0$, there exist a constant $C(\varepsilon,k)$
depending (in addition) on $M$, $C_{min}$ and the dimension only, and $t_\varepsilon>0$ such that
$$
\ent{\mu}{\P_{kt} \nu} \leq \frac{C(\varepsilon,k)}{\log^{k(1-\varepsilon)}(1/\xi(t))} \, ,
$$
for all $t>t_\varepsilon$.
\end{ethm}

Before proving the theorem we need a preliminary result.
Recall first that for all non-negative functions $f,g$ we have
 $\ent{\mu}{f+g} \leq \ent{\mu}{f} + \ent{\mu}{g}$. Then for $h\geq 0$,
applying this with $f=\Pt (h \1_{h\leq K})$ and $g=\Pt (h \1_{h>K})$, and using the fact that
entropy is decaying along the semi-group, we obtain that
\begin{equation}\label{eq-entdecom}
\ent{\mu}{\Pt h} \, \leq \, \ent{\mu}{\Pt (h \1_{h\leq K})}+\ent{\mu}{h \1_{h> K}} \, ,
\end{equation}
for all $K>0$. The next Lemma explains how control the second term of the right hand side
of~\eqref{eq-entdecom} using the estimate of the Proposition~\eqref{prop-royer}.

\begin{elem}\label{lem-entdecay}
Let $h$ be a density of probability with respect to $\mu$.
Assume that there exists $c>0$ such that for all $p>1$,
$$
\left(\int h \, \log_+^p h \, d\mu\right)^{\frac 1p} \leq c p.
$$
For $K\geq e^2$, if $\ent{\mu}{h}\leq \frac {1}{2e} \, \log K$ then we get
\begin{equation}\label{eq-ent1} \ent{\mu}{h \1_{h> K}} \, \leq \, (ec+2) \,
\frac{\ent{\mu}{h}}{\log K} \, \log\left(\frac{\log K}{\ent{\mu}{h}}\right) \, .
\end{equation}
\end{elem}
\begin{eproof}
 It is easily seen (see e.g. \cite[Lemma 3.4]{gui-ca}) that if $K\geq e^2$,
\begin{equation}
\label{eq-cg}
\int \1_{h>K} h d\mu \leq \frac{2}{\log K} \, \ent{\mu}{h} \, .
\end{equation}
 Hence
\begin{eqnarray}\label{eq-entcontrol}
\int \, h \log h \, \1_{h>K} d\mu & \leq & \left(\int \, h  \, \1_{h>K} d\mu\right)^{\frac{p-1}{p}}
\, \left(\int \, h \log_+^p(h) \,  d\mu\right)^{\frac 1p} \\ & \leq & c \, p
\,\left(\frac{\ent{\mu}{h}}{\log K}\right)^{\frac{p-1}{p}}   \leq  ce \,
\frac{\ent{\mu}{h}}{\log K} \, \log\left(\frac{\log K}{\ent{\mu}{h}}\right) \nonumber
\end{eqnarray}
provided $\ent{\mu}{h}\leq \frac 1e \, \log K$. The last inequality is obtained by an optimization
upon $p$ (for which we need $\ent{\mu}{h}\leq \frac 1e \, \log K$).

If $\ent{\mu}{h}\leq \frac {1}{2e} \, \log K$,
$$
- \, \left(\int \1_{h>K} h d\mu\right)
\log\left(\int \1_{h>K} h d\mu\right) \leq \, - \, \left(\frac{2}{\log K} \, \ent{\mu}{h}\right)
\log \left(\frac{2}{\log K} \, \ent{\mu}{h}\right),
$$
using \eqref{eq-cg}, so that we have finished the proof.
\end{eproof}

\emph{\textbf{Proof of Theorem~\ref{thm-entdecay}}}\\\proofbegin~
Let $h=\Ps \nu$. According to \eqref{eq-entdecom}, Proposition \ref{prop-decay1} and Lemma
\ref{lem-entdecay}, it holds for all $t>s>0$,
$$
\ent{\mu}{\Pt \nu} \, \leq \, K \xi(t-s) \, + \, c_s \,
\frac{H}{\log K} \, \log\left(\frac{\log K}{H}\right) \, ,
$$
 where $H=\ent{\mu}{h}$, provided $K$
is large enough. Since $H$ can be bounded from above by a quantity $H_0$ depending on $M$,
$C_{min}$ and the dimension only, we may choose $K>K_1$ independent of $H$.

Choosing $K \, = \, c \, \frac{H_0}{\xi(t-s)} \,
\frac{1}{1+\log_+\left(\frac{H}{\xi(t-s)}\right)}$, we obtain
\begin{equation}\label{eq-dec1}
\ent{\mu}{\Pt \nu} \leq C \, \frac{1+\log_+\left(\log_+(1/\xi(t-s))\right)}{1+\log_+(1/\xi(t-s))}
\,  .
\end{equation}
It follows that, for all $1\geq \varepsilon >0$ there exists some $t_{\varepsilon}$ such that for
$t\geq t_\varepsilon$
\begin{equation}\label{eq-dec2}
\ent{\mu}{\Pt \nu} \leq \frac{C}{\log^{1-\varepsilon}(1/\xi(t))} \, .
\end{equation}
Using again \eqref{eq-entdecom} and \eqref{eq-ent1} (we may choose $c=c_s$ for all $t\geq s$) we
may write
\begin{eqnarray*}
\ent{\mu}{\P_{2t} \nu} & \leq & K \xi(t) \, + \, c \, \frac{\ent{\mu}{\Pt \nu}}{\log K} \,
\log\left(\frac{\log K}{\ent{\mu}{\Pt \nu}}\right) \\ & \leq & K \xi(t) \, + \, \frac{c c'}{\log K
\, \log^{1-2\varepsilon}(1/\xi(t))} \, + \, \frac{c \, \log \log_+K}{\log K \,
\log^{1-\varepsilon}(1/\xi(t))}
\end{eqnarray*}
where we have used $y \log (1/y) \leq c' y^{1-\varepsilon}$ for $y\leq 1/e$. Hence choosing
$K=1/\xi(t) \log^2(1/\xi(t))$ we obtain a bound like
$$
\ent{\mu}{\P_{2t} \nu} \leq\frac{C}{\log^{2-2\varepsilon}(1/\xi(t))} \, ,
$$
for $t$ large enough. Note that $C$ depend on $\varepsilon$. We may iterate the method
and get the result.
\proofend

Of course this result is not totally satisfactory, but it indicates that the decay of entropy is
faster than any $1/\log^{k(1-\varepsilon)}(1/\xi(t/k))$.

\begin{eex}
Let us study the two classical examples we already mentioned. To be rigorous $|t|:=\sqrt{1+t^2}$ in
what follows (to ensure the required regularity), so that \eqref{condition} is satisfied.
\begin{itemize}
\item For $\alpha>0$, the measure $dm_\alpha(t)= Z_\alpha(1+|t|)^{-1-\alpha}dt,\, t\in \dR$
satisfies the weak logarithmic Sobolev inequality with
$$
\forall s\in(0,1),\qquad \beta_{WL}(s)=C\frac{\PAR{\log 1/s}^{1+2/\al}}{s^{2/\al}},
$$
for some constant $C>0$. Hence,

$$
\xi(t) = \frac{c_{\alpha}}{t^{\alpha/2} \, \log^{1+\alpha}(t)}
$$
for large $t$, and
$$
\ent{m_\alpha}{\P_{kt} \nu} \leq  \,\frac{C_{\alpha,k,\e}}{\log^{k(1-\varepsilon)}(t)} \, .
$$
 Notice that, if roughly
the rate of decay does not depend on $\alpha$ (it is faster than any $\log^k(t)$), the dependence
on $\alpha$ of all constants shows that this regime is attained for smaller $t$ when $\alpha$
increases.

\item For $\al\in(0,2)$, the measure $d\mu_\al(t)=Z_\al e^{-|t|^{\al}}dt,\; t\in \dR$, ($Z_\al$ is
a normalization constant) satisfies the weak logarithmic Sobolev inequality  with
$\beta_{WL}(s)=C\PAR{\log 1/s}^{(2-\al)/\al}$, $C>0$. Hence $\xi(t)=c \, e^{- d t^{\alpha/2}}$ and
for $t$ large enough,
$$
\ent{\mu_\alpha}{\P_{kt} \nu} \leq \frac{C_{\alpha,k}}{1+t^{(\alpha/2)(k-\varepsilon)}} \, .
$$ Of course this
result is not satisfactory for $\alpha\geq 2$ where we know that the decay is exponential. See
below for an improvement.
\end{itemize}
\end{eex}

If we replace Proposition \ref{prop-decay1} or Proposition~\ref{thm-resls} we can greatly improve
the previous results. Let us describe the latter situation.

\begin{ethm}\label{thm-decaylatala}
In the situation of Example~\ref{cor-rlslo} (i.e. the Latala-Oleszkiewicz situation) and Theorem
\ref{thm-entdecay}, there exists $s>0$ such that for all $1\geq \varepsilon>0 $ one can find
$T_\varepsilon$ in such a way that for $t\geq T_\varepsilon$,
 $$
 \ent{\mu}{\P_{t+s} \nu} \leq \, e^{1- t^{\frac{(1-\varepsilon)\alpha}{2-\varepsilon \alpha}}} \, .
 $$
 In particular for $\alpha=2$
relative entropy is exponentially decaying.
\end{ethm}

\begin{eproof}
The beginning of the proof is similar to the one of Theorem \ref{thm-entdecay} but replacing the
estimate of Proposition \ref{prop-decay1} by the one of Example~\ref{cor-rlslo} (in particular we
may take $K=+\infty$ if $\alpha=2$). The first step yields $$H_t \, := \, \ent{\mu}{\P_{t+s} \nu}
\leq \, \frac{C (1+\log_+^{\frac{\alpha}{2-\alpha}}(t))}{1+t^{\frac{\alpha}{2-\alpha}}} \, H
\log(1/H) \, .$$ Let us choose $s$ in such a way that $H\leq 1/e$, i.e. $H \log(1/H)\leq 1$. Then
$$H_{2t} \, \leq \, \frac{C (1+\log_+^{\frac{\alpha}{2-\alpha}}(t))}{1+t^{\frac{\alpha}{2-\alpha}}}
\, H_t \log(1/H_t) \, \leq \, \left(\frac{C
(1+\log_+^{\frac{\alpha}{2-\alpha}}(t))}{1+t^{\frac{\alpha}{2-\alpha}}}\right)^2 \,
\log(1+t^{\frac{\alpha}{2-\alpha}}) \, ,$$ provided $C\geq 1$ that we can assume. Iterating the
procedure we get
\begin{eqnarray*}
H_{kt} & \leq & \left(\frac{C
(1+\log_+^{\frac{\alpha}{2-\alpha}}(t))}{1+t^{\frac{\alpha}{2-\alpha}}}\right)^k \,
\prod_{j=1}^{k-1} \, \log\left((1+t^{\frac{\alpha}{2-\alpha}})^j\right) \\ & \leq & \left(\frac{C
\,  (1+\log_+^{\frac{\alpha}{2-\alpha}}(t)) \,
\log(1+t^{\frac{\alpha}{2-\alpha}})}{1+t^{\frac{\alpha}{2-\alpha}}}\right)^k \,
\frac{(k-1)!}{\log(1+t^{\frac{\alpha}{2-\alpha}})}
\end{eqnarray*}
Now, we may find $t_\varepsilon$ such that for $t\geq t_\varepsilon$, $$\frac{C \,
(1+\log_+^{\frac{\alpha}{2-\alpha}}(t)) \,
\log(1+t^{\frac{\alpha}{2-\alpha}})}{1+t^{\frac{\alpha}{2-\alpha}}} \, \leq \,
\frac{1}{t^{\frac{\alpha}{2-\alpha} (1-\varepsilon)}} \, ,$$ and
$\log(1+t^{\frac{\alpha}{2-\alpha}}) \geq 1$, so that $$H_{kt} \, \leq \, \left(\frac{k}{e \,
t^{\frac{\alpha}{2-\alpha} (1-\varepsilon)}}\right)^k$$ as soon as $k$ is large enough (for
$(k-1)!\leq (k/e)^k$). Choosing $t=k^{(2-\alpha)/\alpha (1-\varepsilon)}$ (hence $k$ large enough
for $t$ to be greater than $t_{\varepsilon}$) we obtain that $H_u\leq e^{-k}$ for
$u=k^{\frac{2-\varepsilon \alpha}{(1-\varepsilon)\alpha}}$, i.e. $H_t \, \leq \, e \, e^{-
t^{\frac{(1-\varepsilon)\alpha}{2-\varepsilon \alpha}}}$.
\end{eproof}

Of course the statement of the Theorem is not sharp (we have bounded some logarithm by some power)
but it is tractable and shows that (up to some $\varepsilon$) the decay is similar to $\xi$. Of
course we are able to derive a similar (but not very explicit) result with the general bound ($A$)
in Proposition~\ref{thm-resls}.
\medskip

\subsection{Comparison results and convergence in total variation distance}

It is interesting to see what can be done by using the usual Poincar\'e inequality. Indeed recall
that $\ent{\mu}{g} \leq \varf{\mu}(g)/\int g d\mu$ for a nonnegative $g$. Using this with $g=\Pt (h
\1_{h\leq K})$, using also \eqref{eq-entdecom} and Poincar\'e yield a decay
$$\ent{\mu}{\Pt \nu} \leq C \frac{1+\log_+(t)}{1+t}$$ that is a slightly better result than the one
we may obtain at the first step of the previous method (up to a $\log_+(t)$ factor) in this
situation (corresponding to $\alpha=1$). But iterating the procedure also yields a polynomial
decay. Nevertheless if $\P_s \nu \in \dL^2(\mu)$ for some $s$, we obtain an exponential decay. It
is thus particularly interesting to study stronger integrability condition.

It turns out that Royer's method furnishes a much better result (in a sense) than the one shown in
Proposition \ref{prop-royer}, namely
\begin{eprop}\label{prop-royer2}
Under the hypotheses of Proposition \ref{prop-royer}, for all $t>0$, all $x\in \R^n$, $$\int (\P_t
\delta_x)^2 \, d\mu \, \leq \, (2 \pi t)^{- \frac n2} \, e^{C_{min} t} \, e^{2 V(x)} \, .$$
\end{eprop}

This result can be shown exactly as Proposition \ref{prop-royer} replacing the convex function $u
\mapsto u \log_+^p(u)$ by $u \mapsto u^2$ (see e.g. \cite{gui-ca2}). It shows that for an initial
condition $\nu$, a sufficient condition for $\P_t \nu \in \L^2(\mu)$ (for $t>0$) is $$\int e^{2V}
\, d\nu \, < \, + \infty \, .$$ This is of course a very strong assumption. In particular, it has
been shown by P.A. Zitt (\cite{Zitt}), that Theorem \ref{thm-entdecay} can be used to show the
absence of phase transitions in some infinite dimensional situations, while the control in the
previous Proposition is not useful.

We shall study an example in the next subsection, showing that actually, one can expect a still
better integrability for $\P_t \delta_x$.
\bigskip

To finish this section we shall now discuss the weaker convergence in total variation distance.

Denoting again $h=\Ps \nu$, we thus have for $K>0$
\begin{eqnarray}\label{eq-TV}
\int  |\Pt h - 1| d\mu & \leq & \int  |\Pt (h\wedge K) - \Pt h| d\mu + \int  |\Pt (h\wedge K) -
\int (h\wedge K) d\mu| d\mu + |\int (h\wedge K) d\mu - 1| \nonumber \\ & \leq & \int  |\Pt (h\wedge
K) - \int (h\wedge K) d\mu| d\mu + 2 \int (h-K) \1_{h\geq K} d\mu
\end{eqnarray}
where we have used the fact that $\Pt$ is a contraction in $L^1$.
The second term in the right hand sum is going to 0 when $K$ goes to $+\infty$, while the first
term can be controlled either by $\sqrt{\varf{\mu}(\Pt (h\wedge K))}$ or by $\sqrt{2 (\int (h\wedge
K) d\mu) \, \ent{\mu}{\Pt (h\wedge K)}}$ according respectively to Cauchy-Schwarz and to Pinsker
inequality. In both cases, {\bf WPI} or {\bf WLSI} inequalities imply that $\Pt \nu$ goes to $\mu$
in total variation distance, for all initial $\nu$.
\medskip

If we want a rate of convergence, we immediately see that {\bf WPI} will furnish a better rate
than {\bf WLSI} for the $\mu$ that do not satisfy Poincar\'e inequality. If $\mu$ satisfies a
Poincar\'e inequality with constant $C_P$ then $$\varf{\mu}(\Pt (h\wedge K)) \leq K e^{-t/C_P} \,
,$$ so that the optimal $K$ is  given (up to a factor 2) by $2 \int (h-K) \1_{h\geq K} d\mu =
K^{\frac 12} e^{-t/2C_P}$. In particular if \eqref{condition} holds, $$2 \int (h-K) \1_{h\geq K}
d\mu \leq \frac{2 C(p)}{\log^p(K)}$$ for $K>1$ and $p\geq 1$, so that we obtain $\parallel
\P_{t+s} \nu - \mu\parallel_{TV} \, \leq \, \kappa(p)/t^p \, $ for all $s>0$, $p\geq 1$, where
$\kappa$ depends on $s$, $C_{min}$, $p$, $M$ and the dimension. But if we directly use Theorem
\ref{thm-decaylatala} and Pinsker we have the much better $\parallel \P_{t+s} \nu -
\mu\parallel_{TV} \, \leq \, \kappa \, e^{- \, \frac 12 \,
t^{\frac{(1-\varepsilon)\alpha}{2-\varepsilon \alpha}}}$ at least for $s$ large enough. In
particular for $\alpha=1$ we obtain a faster decay. Once again, if $\| \P_{s}\nu\|_\infty$ is
finite for some positive $s$ then one should use the entropic convergence of
Proposition~\ref{thm-resls} to get an exponential decay.
\medskip

\subsection{Example(s)}

In the previous subsections, we have seen that finite entropy conditions are quite natural for the
law of the diffusion at any positive time, but that starting from an initial Dirac mass, we
immediately reach $\L^2(\mu)$. Before to study examples indicating that one can expect much
better, we shall give a generic example showing that some natural measures $\nu$ never satisfy
$\Ps \nu \in \dL^2(\mu)$, but satisfy the conditions in Proposition \ref{prop-royer}.
\medskip

Consider $V$ such that for all $\lambda>0$ , $\int e^{- \lambda V} dx < +\infty$. Let
$d\mu=e^{-2V}dx$ and $d\nu=e^{-(2-\varepsilon)V}/Z_\varepsilon \, dx$ so that $d\nu/d\mu:=h =
Z_\varepsilon \, e^{\varepsilon V} \notin \dL^2(\mu)$ for $2>\varepsilon >1$, but $\int
e^{\frac{2-\varepsilon}{2} \, V} d\nu <+\infty$. Set $G=e^V= h^{\frac{1}{\varepsilon}}$.

If $\Ps h \in \dL^2(\mu)$ for some $s>0$, then $\Ps G \in \dL^{2 \varepsilon}(\mu)$. If
\eqref{condition} holds, it follows from \cite[Theorem 2.8]{cat5} that $\mu$ satisfies a
logarithmic Sobolev inequality. Thus if it is not the case, $\Ps h \notin \dL^2(\mu)$ for all $s\geq 0$,
while if \eqref{eq-initial} is satisfied (for instance for $V(y)=|y|^{\alpha}$, $1\leq \alpha < 2$
see below) $\nu$ satisfies the conditions in Proposition \ref{prop-royer}.

This example shows that the set of initial measures satisfying the conditions in the previous
subsection but not the necessary conditions to simply apply Poincar\'e is non empty.
\medskip

We shall go further, and for simplicity we shall only consider the measures $\mu_\alpha$ for
$\alpha \geq 1$, and essentially discuss the case $\alpha=1$.

First of all notice that if $1\leq \alpha \leq 2$, $$|y|^\alpha \leq 2^{\alpha - 1} (|x|^\alpha +
|y-x|^2 + 1)$$ so that \eqref{eq-initial} is satisfied. Hence as soon as $\int e^{\lambda
|x|^\alpha} \nu(dx) < +\infty$ for some $\lambda >0$, we may apply all the results of the previous
subsection. We shall now give a precise description of $h = \Ps \delta_x$. This will allow us to
give a similar sufficient condition for $\Ps \nu$ to belong to $\dL^2(\mu)$.
\medskip

We thus consider (in one dimension)
\begin{equation}\label{eq-expo}
dX_t \, = \, dB_t \, -  \, \textrm{sign}(X_t)dt \quad , \quad X_0=x \, ,
\end{equation}
corresponding to $\alpha=1$. Elementary stochastic calculus (inspired by the first sections of
\cite{GHR}) furnishes
\begin{eqnarray*}
\dE [f(X_t)] & = & \dE \left[f(x+B_t) \, e^{-\frac t2}\, \exp \left(-\int_0^t \textrm{sign}(x+B_s)
dB_s\right)\right] \\ & = & e^{|x|} \, e^{- \frac t2} \, \dE \left[f(x-W_t) \, \exp \left(-|W_t -
x| + L_t^x\right) \right]
\end{eqnarray*}
where $W_s=-B_s$ is a new Brownian motion with local time at $x$ denoted by $L_s^x$. Now as usual
we introduce the hitting time of $x$ of $(W_s)$ denoted by $T_x$, and the supremum $S_t=\sup_{0\leq s \leq t}
W_s$. We also assume here that $x > 0$. Then
\begin{eqnarray*}
\dE [f(X_t)] & = & \dE [f(X_t) \, \1_{t\leq T_x}] + \dE [f(X_t) \, \1_{t> T_x}] \\ & = & e^{|x|} \,
e^{- \frac t2} \, \dE [f(x-W_t) \, \1_{S_t\leq x} \, e^{W_t-x}] \, + \, e^{-\frac t2} \, \dE
[\1_{S_t > x} \dE[f(B'_{t-T_x}) \, \exp \left(- |B'_{t-T_x}|+L'_{t-T_x}\right)]]
\end{eqnarray*}
where $B'$ is a Brownian motion independent of $W$ and $L'$ its local time at 0.

For the first term, we know that the joint law of $(W_t,S_t)$ is given by the density $$(w,s)
\mapsto \1_{w\leq s} \, \sqrt{2/\pi t^3} \, (2s-w) \, \exp (- (2s-w)^2/2t) $$ so that (recall $x>0$) $$\dE [f(X_t)
\, \1_{t\leq T_x}] = \int f(u) \left( \1_{u\geq 0} \sqrt{2/\pi t} \, e^{-\frac t2} \, e^x \, e^{-u}
\, \left(e^{-(x-u)^2/2t} - e^{-(x+u)^2/2t}\right) \right) du \, .$$

For the second term, we know that the law of $T_x$ is given by the density $$T \mapsto x \,
\sqrt{1/2\pi T^3} \, e^{-x^2/2T}$$ and that $(|B'_s|,L'_s)$ has the same law as $(S'_s-B'_s,S'_s)$
so that (noting that only the even part of $f$ has to be considered)
$$\dE [f(X_t) \, \1_{t> T_x}]= e^{- \frac t2} \, \iiint \1_{0 <T < t} \1_{u>0} \1_{v>u} \left(\frac{f(u)+f(-u)}{2}\right)
 g(T,u,v) \, du dv dT \, ,$$
with $$g(T,u,v)=\sqrt{1/2\pi T^3} \, \sqrt{2/\pi (t-T)^3} \,  v \, e^v \, e^{-2u} \,
e^{-v^2/2(t-T)} \, e^{-x^2/2T} \, .$$ But $$Q:= \int_0^t \int_u^{+\infty} \, \sqrt{1/2\pi T^3} \,
\sqrt{2/\pi (t-T)^3} \,  v \, e^v \, e^{-v^2/2(t-T)} \, e^{-x^2/2T} \, dv dT$$ is such that
\begin{eqnarray*}
Q & \leq & \int_0^t \, \sqrt{1/2\pi T^3} \, \left(\sqrt{2/\pi (t-T)}\,  e^{u} \, e^{-u^2/2(t-T)}+ 2
e^{t-T}\right) \, e^{-x^2/2T} \,  dT \\ & \leq & \int_0^t \, \sqrt{1/2\pi T^3} \, \left(\sqrt{2/\pi
(t-T)}\,  e^{t/2}+ 2 e^{t-T}\right) \, e^{-x^2/2T} \,  dT \\ & \leq & C(t)
\end{eqnarray*}
independently of $x$. The first inequality is obtained by performing an integration by parts in
$v$, the second one by bounding $e^{u} \, e^{-u^2/2(t-T)}$ and the final one by bounding separately
$\int_0^{t/2}$ and $\int_{t/2}^t$. We thus see that $$\dE [f(X_t) \, \1_{t> T_x}]=C'(t) \, \int
f(u) \, e^{-2|u|} \, g(u) \, du$$ where $g$ is bounded.

Putting all this together we have obtained the following
\begin{equation}\label{eq-densite}
(\Pt \delta_x)(u) = c(t) \left( \1_{u\geq 0} \, e^x \, e^{u} \, \left(e^{-(x-u)^2/2t} -
e^{-(x+u)^2/2t}\right)\right) + C'(t) g(u)
\end{equation}
for all $x>0$. A similar result holds for $x<0$, while $\Pt \delta_0$ is bounded. Of course the
previous \eqref{eq-densite} shows that for a fixed $x$, $\Pt \delta_x$ is bounded. This result is
not so surprising. Indeed for $\alpha=2$ (more precisely for the normalized gaussian measure i.e.
the Ornstein-Uhlenbeck process) $(\Pt \delta_x)(u)= c(t) \, e^{(1-e^{-t})x^2/2(1-e^{-t})} \,  e^{-
(e^{-t/2}u-x)^2/2(1-e^{-t})}$ is bounded too. One may adapt our proof and Proposition 4 in
\cite{GHR} in order to show that a similar result actually holds for all $1\leq \alpha\leq 2$.

But \eqref{eq-densite} allows us to look at more general $\Pt \nu$. In particular we see that $\Pt
\nu \in \dL^2(\mu)$ if and only if \begin{equation}\label{nosg}\int_{u>0} \, \left(\int_{x>0} \, e^x \, e^{-(u-x)^2/2t} \,
\nu(dx)\right)^2 \, du \, < \, +\infty\end{equation} and a similar property is available on the negative real
numbers. We then easily recover and complete the discussion at the beginning of this subsection,
i.e. if $d\nu=e^{-\lambda |x|} dx/Z$ , $\Pt \nu \notin \dL^2(\mu)$ if $\lambda \leq 1$, but belongs
to $\dL^2(\mu)$ if $\lambda>1$.

Let us finally give some discussion concerning the obtainable rate of entropic convergence
depending on the initial measure:
\begin{enumerate}
\item if $\nu=\delta_x$, then $\|\P_{t_0}\delta_x\|_\infty<\infty$ and using respectively
Proposition \ref{prop-decay1}, Proposition~\ref{thm-resls} or Poincar\'e inequality, one gets
$$\ent{\mu}{\P_{t+t_0}\nu}\le C \min\left(  e^{-a\sqrt{t}}\|\P_{t_0}\delta_x\|_\infty, e^{-bt/
(1+\log\|\P_{t_0}\delta_x\|_\infty)},e^{-ct}\|\P_{t_0}\delta_x\|_\infty\right),$$ (note that it
easily extends to the case where $\nu$ has compact support.) \item if $\nu$ does not satisfy
(\ref{nosg}) but for some positive $\lambda$, $\int e^{\lambda |x|}d\nu$  is finite then we can
only use Theorem \ref{thm-decaylatala} to get that for all $\varepsilon>0$, there exists
$T_\varepsilon$ such that for all $t\ge T_\varepsilon$ we have $$\ent{\mu}{\P_{t}\nu}\le
e^{1-t^{1-\varepsilon\over 2-\varepsilon}}.$$

\end{enumerate}

\section{Classical properties of WLSI}
\label{sec-pw}
\subsection{Tensorization}
Let us begin by the following naive procedure of tensorization.

\begin{eprop}
Assume that $\mu$ satisfies a {\bf WLSI} with function $\beta$ and let $n\ge1$. Then the measure
$\mu^n$ satisfies a {\bf WLSI} with function $\beta\Big(\frac{s}{n}\Big)$, for $s>0$.
\end{eprop}

\begin{eproof}
By the sub-additivity property of the entropy we get
 $$ \ent{\mu^{n}}{f} \leq \sum_{i=1 }^{n} \int \ent{\mu}{
  f(x_1,\ldots, x_{i-1},\cdot,x_{i+1},\ldots,x_n)} \prod_{j\neq i}d\mu(x_j).$$
 For each $i$ we get for all $\PAR{x_1,\ldots, x_{i-1},x_{i+1},\ldots,x_n}\in M^{n-1}$
 \begin{multline*}
\ent{\mu}{f(x_1,\ldots, x_{i-1},\cdot,x_{i+1},\ldots,x_n)}\leq\\
\beta(s) \int |\nabla_i f|^{2}(x_1,\ldots,y_i,\ldots,x_n) d\mu(y_i)+s \,
  \osc{ f(x_1,\ldots,\cdot,\ldots,x_n)}^{2},
\end{multline*}
 It yields
 $\forall s>0,\qquad \ent{\mu^{n}}{f} \le \beta(s) \int |\nabla f|^{2}d\mu^{n}+ns \, {\bf Osc}^2(f).$
  \end{eproof}

The tensorization result above is of course the same as the one in \cite{ca-ba-ro2} for weak
Poincar\'e inequality. As explained in Section 5 of this paper, one cannot expect a better result
beyond the exponential case. However as we have already seen, {\bf WLSI} may take place between
the exponential and the gaussian regime (when {\bf GBI} holds), so that we obtain this corollary:
\begin{ecor}
If $\mu_i$ ($1\leq i\leq n$) satisfy a {\bf WLSI} with the same function $\beta_{WL}$ satisfying
the hypotheses in Proposition \ref{prop-becwl}, then the tensor product $\otimes_{i=1}^{n} \mu_i$
satisfies a {\bf WLSI} with function
$$
\beta^n_{WL}(u)=C \, \beta_{WL}(C'u )
$$
where $C,C'$ are constants which don't depend on $n$.
\end{ecor}
\begin{eproof}
It is enough to use both parts of Proposition \ref{prop-becwl} and the (exact) tensorization
property of {\bf GBI}. One can see~\cite{LO00} for the proof of the tensorization of {\bf GBI}.
\end{eproof}

Among the most important consequences of functional inequalities, one find concentration of
measure and isoperimetric profile. Unfortunately weak inequalities are not easily tractable to
derive results in this direction (due to the Oscillation term). However results for {\bf WPI} are
contained in \cite{r-w,ca-ba-ro2} with a particular interest in dimension dependence in the
latter. Actually we do not succeed in deriving similar estimates starting from {\bf WLSI}, as
Herbst's argument or Aida-Masuda-Shigekawa iteration argument are more intricate and we can only
recover weak Poincar\'e non optimal concentration rate.

The situation is still worse (from the {\bf WLSI} point of view) when a {\bf SPI} holds. In this
case various (more or less explicit) results have been obtained. Let us mention on one hand
\cite{w1} Section~6, \cite{GWa02} Section 5 (using super Poincar\'e) and \cite{w2} Corollary 2.4
(using {\bf GBI}), on the other hand \cite{ca-ba-ro} Section~6 (using {\bf GBI}) and Section 8
(using $F$-Sobolev inequalities) and \cite{ca-ba-ro3} Theorem 12 for an improvement of \cite{w1}
Section~6. The previous result may be used in conjunction with the above mentioned results to get
dimension free concentration (or isoperimetric) results, completing thus the transportation
approach presented before.

\subsection{Perturbation}

Among the methods used to obtain functional inequalities, an efficient one is to perturb measures
satisfying themselves some functional inequalities. The most known result in this direction was
first obtained by Holley and Stroock who showed that a logarithmic Sobolev inequality is stable under a
log-bounded perturbation. The same is true for a {\bf SPI} (using the
related {\bf GBI} \cite[Proposition 2.5]{w2}), and actually one can replace the
bounded assumption by a Lipschitz assumption (this was shown by Miclo for logarithmic Sobolev, and by Wang
\cite[Proposition 2.6]{w2} for a {\bf SPI}).

For the {\bf WPI}, a similar result is shown in \cite[Theorem 6.1]{r-w}. Actually
this result shows that one can consider non bounded perturbation, but with very strong
integrability assumptions, the final result being far to be explicit. For {\bf WLSI} we may state
\begin{eprop}
\label{prop-bounded}
Suppose that $\mu$ satisfies a {\bf WLSI} with function $\beta_{WL}$. Let
$\nu_V={e^{V}} \mu/Z_V,$ where $Z_V=\int e^Vd\mu$ and assume that $V$ is bounded on $M$.

Then $\nu_V$
satisfies a {\bf WLSI} with function
$$
\beta_{WL}^V(u)= e^{2{\bf Osc}(V)} \, \beta_{WL}(u e^{- {\bf Osc}(V)}) \, .
$$
We may replace {\bf WLSI} by {\bf WPI} replacing $\beta_{WL}$ by $\beta_{WP}$, or by
{\bf SPI} with
$$
\beta_{SP}^V(u)=e^{2{\bf Osc}(V)} \, \beta_{SP}(u e^{-2 {\bf Osc}(V)}).
$$
\end{eprop}

\begin{eproof}
Recall that $\ent{\nu_V}{f^2} \, \leq \, e^{{\bf Osc}(V)} \, \ent{\mu}{f^2}$.  Applying {\bf WLSI} for
$\mu$ yields
\begin{eqnarray*}
\ent{\nu_V}{f^2} & \leq & e^{{\bf Osc}(V)} \, \left(\beta_{WL}(s)\int\ABS{\nabla f}^2d\mu+s \, {\bf
Osc}^2(f)\right) \\ & \leq & e^{2{\bf Osc}(V)} \, \beta_{WL}\PAR{u e^{- {\bf Osc}(V)}} \, \int\ABS{\nabla
f}^2d\nu_V + u \, {\bf Osc}^2(f) \, ,
\end{eqnarray*}
which is exactly the first statement. The second one is similar since
$\varf{\nu_V}(f) \leq
e^{{\bf Osc}(V)} \varf{\mu}(f).$
 For {\bf SPI} the proof is immediate.
\end{eproof}

The second way to get perturbation results is to use a natural isometry between $\dL^2$ spaces. For
notational convenience we assume now that $\nu_V = e^{-2V} \mu$. Then $g \mapsto f:=e^{-V} g$ is an
isometry between $\dL^2(\nu_V)$ and $\dL^2(\mu)$. It is thus immediate that on one hand

\begin{equation}\label{eq-perturb}
\ent{\nu_V}{g^2} \, = \, \ent{\mu}{f^2} \, + \, 2 \int g^2 V \, d\nu_V \, .
\end{equation}

On the other hand, an integration by parts yields
\begin{equation}\label{eq-perturbgrad}
\int\ABS{\nabla f}^2 d\mu = \int\ABS{\nabla g}^2 d\nu_V \, + \, \int g^2 \, \left(2 LV -
\ABS{\nabla V}^2\right) d\nu_V \, ,
\end{equation}
where $L$ is the generator of $P_t$ reversible for $\mu$.

Combining these two facts, yields perturbation results for logarithmic Sobolev inequalities (the idea goes
back to Rosen \cite{Ros}, and was used in \cite{Carl,cat5}). In order to see how to use it in our
framework, we shall first introduce some notation.

\begin{edefi}
\label{def-Gset} Let $G$ be a positive continuous function defined on $\dR^+$. We shall say that a
smooth $V$ is ($G,\mu$)-good, if $V(x) \to +\infty$ as $|x| \to +\infty$ and if there exists
$A\geq 0$ such that one has for any $x$ such that ${V(x)\geq A}$,
$$
\ABS{\nabla V}^2(x) - 2 LV(x)\geq G(V(x)) \, .
$$
\end{edefi}
Our first general result is a bounded (but not log-bounded) perturbation result.

\begin{eprop}\label{prop-Wit}
Let $\mu$ be a positive measure (not a necessarily probability measure) satisfying a {\bf WLSI}
with continuous function $\beta_{WL}$. Let $V$ be ($G,\mu$)-good, such that $\nu_V=e^{-2V}\mu$ is
a probability measure.
\smallskip

Then for all $u>0$ and $b\geq A$ the following inequality holds for any $g\in H^1(E,\mu)$,
$$
\ent{\nu_V}{g^2} \, \leq \, C(u,b)  \, \int\ABS{\nabla g}^2 d\nu_V \, + \, D(u,b) \, {\bf Osc}^2(g) \, ,
$$
 with
\begin{equation}\label{eq-C}
C(u,b) \, = \, h(b) +\left(2+2A+M(V)h(b)\right) \, \beta_{WP}^V(u) \, ,
\end{equation}
\begin{equation}\label{eq-D}
D(u,b) \, = \, s_b \, e^{- 2 \inf V} + \left(2+2A+M(V)h(b))\right) \, u + \int_{\{V\geq b\}} 2V
d\nu_V \, ,
\end{equation}
\smallskip

where $h(b):= \sup_{\{A\leq z \leq b\}} \, \frac{2z}{G(z)}$ , $s_b:=\inf \, \{s>0 \, , \,
\beta_{WL}(s)\leq
 h(b)\}$,
$$
M(V):=\sup_{\{V\leq A\}} \, (2LV - \ABS{\nabla V}^2),
$$
(which is finite) and
$\beta_{WP}^V$ is the best function such that $\nu_V$ satisfies {\bf WPI} (if it does not take
$\beta_{WP}^V(u)=+\infty$ for small $u$).
\end{eprop}

\begin{eproof}
First according to Rothaus inequality, we may assume that $\int g d\nu_V=0$ up to $2
\varf{\nu_V}(g)$.
Applying {\bf WLSI} in \eqref{eq-perturb} and \eqref{eq-perturbgrad} we get for all $s>0$,
\begin{multline}\label{eq-entperturb}
\ent{\nu_V}{g^2} \, \leq \, \beta_{WL}(s) \int\ABS{\nabla g}^2 d\nu_V \, + \\\, \int g^2 \,
\left(\beta_{WL}(s) \left(2 LV - \ABS{\nabla V}^2\right) + 2V\right) d\nu_V \, + \, s {\bf Osc}^2(g
e^{-V}) \, .
\end{multline}
Note that if $\beta_{WL}$ is bounded, we may replace it by any  $\beta (s)\geq \beta_{WL}(0)$.
\begin{itemize}
\item  On $\{V\leq A\}$, the second integrand is bounded by $\left(\beta_{WL}(s)M(V)+2A\right) \,
\varf{\nu_V}(g)$, and can be controlled (together with the term $2\var{\nu_V}{g}$ coming from Rothaus
inequality) with the {\bf WPI} for the measure $\nu_V$.
\item On $\{b\geq V\geq A\}$, we choose $s=s_b$ then  the second integrand is
non-positive.
\item On $\{b\leq V\}$, $2LV - \ABS{\nabla V}^2$ is still non-positive, so that the second
integrand is bounded by $$\int_{\{V\geq b\}} 2V g^2 \, d\nu_V \, \leq \, \left(\int_{\{V\geq b\}}
2V d\nu_V\right) \, {\bf Osc}^2(g) \,,$$
\end{itemize}
since $\int gd\nu_V=0$.
\end{eproof}

For this proposition to be useful, we must choose $u$ and $b$ in such a way that $D(u,b) \to 0$ as
$b \to +\infty$.

If $\mu$ is a probability measure, $\int e^{2V} d\nu_V=1$ so that if $b>1/2$, $$\int_{\{V\geq b\}}
2V d\nu_V \leq \ent{\nu_V}{\1_{V\geq b}} = \nu_V(V\geq b) \, \log\left(\frac{1}{\nu_V(V\geq
b)}\right)\leq b \, e^{-2b}$$ where we used Markov inequality and the fact that $x \log(1/x)$ is
non decreasing on $[0,1/e]$ for the latter.

If $\mu$ is not bounded, we assume in addition that $\int e^{-pV} d\mu = K(p) < +\infty$ for some
$p<2$, so that a similar argument (changing the constants) yields again $$\int_{\{V\geq b\}} 2V
d\nu_V \leq \nu_V(V\geq b) \, (2/2-p) \, \log\left(\frac{K(p)}{\nu_V(V\geq b)}\right) \leq
(2K(p)/(2-p)) \, b \, e^{(p-2)b}$$ if $b\geq (1+\log(K(p))/(2-p)$.

In both cases, defining $\varepsilon$ as the upper bound, one can find constants $a$ and $a'$
(depending on $p$ if necessary) such that $$b=a \log\left(\frac{a'
\log(1/\varepsilon)}{\varepsilon}\right) \, ,$$and the appropriate choice for $u$ is then
$u=\varepsilon/h(b)$, provided $\beta_{WL}(\varepsilon)\leq h(b)$.

Conversely, if $\beta_{WL}(\varepsilon)\geq h(b)$, $s_b$ is greater than $\varepsilon$ (up to
multiplicative constants) and the good choice is then $u=s_b/h(b)$.
\smallskip

If $h(b)\geq C b$ we obtain that $\beta_{WL}^V(s)$ behaves like a function greater than or equal to
(up to some constants) $\log(1/s) \beta_{WP}^V(s/\log(1/s))$ in the first case, $\beta_{WL}(s)
\beta_{WP}^V(s/\beta_{WL}(s))$ in the second case, with $\beta_{WL}(s)$ larger than $\log(1/s)$ in
the latter case. Hence the result is not better (even worse) than \eqref{eq-ddd} in Proposition
\ref{prop-dd}.
\smallskip

If $h(b)/b \to 0$ as $b\to +\infty$ we obtain the same results, but replacing $\log(1/s)$ by
$h(\log(1/s))$, provided $\beta_{WP}^V$ is not bounded (otherwise $\beta_{WL}^V(s)= C h(\log(1/s))$
for some $C$). Hence if $\beta_{WL}(s)\ll \log(1/s)$ we obtain a better result that the one in
Proposition \ref{prop-dd}, namely $\nu_V$ satisfies {\bf WPI} with a function
$$\beta(s)\geq \frac{h(\log(1/s))}{\log(1/s)} \, \beta_{WP}^V(cs)$$ provided this function is
non-increasing. But if there exists $M$ such that $\beta_{WP}^V(c s)\leq M \beta_{WP}^V(s)$, we may
thus choose $\beta \leq (1/2) \beta_{WP}^V$, which leads to a contradiction since $\beta_{WP}^V$ is
assumed to be the best one. We have thus obtained (recall that we leave some constants away in the
previous argument)

\begin{ecor}\label{cor-ts}
Let $\mu$ be a positive measure (not necessarily bounded) satisfying a {\bf WLSI} with continuous
function $\beta_{WL}$. Let $V$ be ($G,\mu$)-good, such that $\nu_V=e^{-2V}\mu$ is a probability
measure. If $\mu$ is not bounded, we assume in addition that there exists $p<2$ such that $\int
e^{-pV} d\mu <+\infty$.

Assume in addition that
\begin{itemize}
\item $h(b):= \sup_{\{A\leq z \leq b\}} \, \frac{2z}{G(z)}$ is such that $h(b)/b \to 0$ as $b\to
+\infty$, \item $\beta_{WL}(s)/ \log(1/s) \, \to 0$ as $s\to 0$ (that is, if $\mu$ is bounded,
$\mu$ satisfies some {\bf SPI} which is stronger than the usual Poincar\'e inequality).
\end{itemize}

Then $\nu_V$ satisfies a Poincar\'e inequality, and a {\bf WLSI} with function $\beta_{WL}^V(s)=a
h(a' \log(1/s))$ for some constants $a$ and $a'$.
\smallskip

In particular if $G(z)\geq cz$ for large $z$, $\nu_V$ satisfies the usual logarithmic Sobolev inequality.
\end{ecor}
\smallskip

The previous result extends part of the results in \cite{cat5} since we do not assume that $\mu$
satisfies a logarithmic Sobolev inequality.

It has to be noticed that the conditions in Corollary \ref{cor-ts} are far to be optimal for
$\nu_V$ to satisfy Poincar\'e inequality. Indeed if $\mu=dx$ on the euclidean space, it is known
that $G(b)\geq k>0$ for large $b$ is sufficient (i.e. $h$ asymptotically linear) (see \cite{cat5}
for a reference). In the general manifold case with $\mu$ the riemannian measure, Wang (\cite{w0}
Theorem 1.1 and Remark 1) has obtained a beautiful sufficient condition, namely $- L\rho(x) \geq k
> 0$ for $\rho(x)$ large, when $\rho$ is the riemannian distance to some point $o$. In the flat
case, this condition reads $|\nabla V|(x)>k>0$ for $|x|$ large. In the one dimensional case, it is
easy to see that this condition is weaker than our $G(b)\geq k>0$ for large $b$. Wang's condition
thus appears as the best general one, though it is not necessary as shown in one dimension by a
potential $V(x)=x+\sin(x)$ for large $x$. But Wang's approach, based on Cheeger inequality and the
control of local Poincar\'e inequality outside large balls, seems difficult to extend to more
general functional inequalities (though it can be used in particular cases, see \cite{r-w} section
3 and \cite{w1}).

\begin{eex}
For $1<\alpha \leq 2$ and $G(u)=u^{2(1-\frac{1}{\alpha})}$ we recover (here $d\mu=dx$) the same
$\beta_{WL}$ as the one corresponding to the measure $\mu_{\alpha}$ studied at the end of section
\ref{sec-wls}. This furnishes a new proof of some results in \cite{ca-ba-ro} section 7.2. For more
general $G$ the result is linked to the perturbation results in \cite{ca-ba-ro3}.
\end{eex}

\bibliographystyle{alpha}
\newcommand{\etalchar}[1]{$^{#1}$}
\def\cprime{$'$}

\bigskip
\noindent P. Cattiaux: Ecole Polytechnique, CMAP, 91128 Palaiseau Cedex France and Universit\'e
Paris X Nanterre, Equipe MODAL'X, UFR SEGMI, 200 avenue de la
R\'epublique, 92001 Nanterre Cedex, France.\\
\noindent
Email: cattiaux@cmapx.polytechnique.fr

\medskip\noindent
I. Gentil: Universit\'e Paris-Dauphine, CEREMADE, UMR CNRS 7534, Place du Mar\'echal De Lattre De Tassigny,
75775 Paris Cedex 16, France.\\
\noindent
Email: gentil@ceremade.dauphine.fr

\medskip\noindent
A. Guillin:  Ecole Centrale Marseille et LATP
UMR CNRS 6632, Centre de Mathematiques et Informatique
Technopôle Château-Gombert,
39, rue F. Joliot Curie,
13453 Marseille Cedex 13, France.\\
\noindent
Email:   guillin@cmi.univ-mrs.fr

\end{document}